\documentclass[twoside]{amsproc}

\usepackage{amscd,amssymb}

\def\to{\mathchoice
{\longrightarrow}
{\rightarrow}
{\rightarrow}
{\rightarrow}}

\def\mapsto{\mathchoice
{\DOTSB\mapstochar\longrightarrow}
{\DOTSB\mapstochar\rightarrow}
{\DOTSB\mapstochar\rightarrow}
{\DOTSB\mapstochar\rightarrow}}

\def\hookrightarrow{\mathchoice
{\DOTSB\lhook\joinrel\relbar\joinrel\rightarrow}
{\DOTSB\lhook\joinrel\rightarrow}
{\DOTSB\lhook\joinrel\rightarrow}
{\DOTSB\lhook\joinrel\rightarrow}}

\numberwithin{equation}{section}

\newtheorem{theorem}{Theorem}[section]
\newtheorem*{conjecture}{Virasoro Conjecture}
\newtheorem{proposition}[theorem]{Proposition}
\newtheorem{lemma}[theorem]{Lemma}

\theoremstyle{definition}

\newtheorem{definition}[theorem]{Definition}

\newcommand{\Z}{\mathbb{Z}}
\newcommand{\C}{\mathbb{C}}
\newcommand{\N}{\mathbb{N}}
\newcommand{\R}{\mathbb{R}}
\newcommand{\Q}{\mathbb{Q}}

\newcommand{\om}{\omega}

\renewcommand{\phi}{\varphi}
\renewcommand{\o}{\otimes}

\DeclareMathOperator{\End}{End}

\newcommand{\p}{\partial}
\newcommand{\Om}{\Omega}
\renewcommand{\*}{\cdot}

\renewcommand{\[}{{[\![}}
\newcommand{\<}{\langle}
\renewcommand{\>}{\rangle}
\newcommand{\Mbar}{\overline{\mathcal{M}}}
\DeclareMathOperator{\ch}{ch}
\DeclareMathOperator{\Tr}{Tr}
\DeclareMathOperator{\Str}{Str}
\newcommand{\bull}{\bullet}

\DeclareMathOperator{\SL}{SL}

\newcommand{\half}{\tfrac12}
\newcommand{\Id}{I}

\newcommand{\g}{\mathfrak{g}}

\newcommand{\CL}{\mathcal{L}}
\newcommand{\CT}{\mathcal{T}}
\newcommand{\CC}{\mathcal{U}}
\newcommand{\DC}{\mathcal{V}}
\newcommand{\UU}{\mathsf{G}}

\newcommand{\Nov}{\Lambda}
\DeclareMathOperator{\ZE}{ZE}
\DeclareMathOperator{\supp}{supp}
\DeclareMathOperator{\ev}{ev}
\newcommand{\virt}{{\textup{virt}}}
\newcommand{\Cbar}{\overline{\mathcal{C}}}
\newcommand{\CO}{\mathcal{O}}
\newcommand{\CV}{\mathcal{V}}
\newcommand{\CE}{\mathcal{E}}
\newcommand{\PP}{\mathsf{P}}
\newcommand{\DD}{\mathsf{D}}
\DeclareMathOperator{\Todd}{Todd}
\DeclareMathOperator{\ad}{ad}
\DeclareMathOperator{\vdim}{vdim}

\newcommand{\QQ}{\mathcal{A}}
\newcommand{\Dil}{\mathcal{D}}
\newcommand{\PHASE}{\mathcal{H}}
\newcommand{\phase}{\mathsf{H}}
\newcommand{\Euler}{\mathcal{E}}

\DeclareMathOperator{\Res}{Res}
\DeclareMathOperator{\Sl}{sl}

\renewcommand{\u}{\mathbf{u}}
\newcommand{\HH}{\mathsf{K}}

\begin{document}

\title{The Virasoro conjecture for Gromov-Witten invariants}

\author{E. Getzler}

\address{Northwestern University, Evanston, IL 60208-2730, USA}

\email{getzler@math.nwu.edu}

\thanks{This work is partially funded by the NSF under grant DMS-9704320.}


\begin{abstract}
The Virasoro conjecture is a conjectured sequence of relations among the
descendent Gromov-Witten invariants of a smooth projective variety in all
genera; the only varieties for which it is known to hold are a point
(Kontsevich) and Calabi-Yau manifolds of dimension at least three. We
review the statement of the conjecture and its proof in genus $0$,
following Eguchi, Hori and Xiong.
\end{abstract}

\maketitle

\section{Introduction}

Now that there exist rigourous constructions of Gromov-Witten invariants of
smooth projective varieties over $\C$ (and more generally, of compact
symplectic mani\-folds), there is growing interest in calculating them and
studying their properties. One of the most intriguing conjectures in the
subject is the Virasoro conjecture of Eguchi, Hori and Xiong \cite{EHX}.

Let $V$ be a smooth projective variety; the Gromov-Witten invariants
of $V$
$$
\<\tau_{k_1}(x_1)\dots\tau_{k_n}(x_n)\>^V_{g,\beta} \in \Q
$$
are parametrized by cohomology classes $x_i\in H^\bull(V,\Q)$, natural
numbers $k_i$, a genus $g\ge0$ and a degree $\beta\in H_2(V,\Z)$. These
invariants are multilinear in the cohomology classes $x_i$ and graded
symmetric under simultaneous permutation of $x_i$ and $k_i$. We recall the
definition of the Gromov-Witten invariants in Section~\ref{definition}.

Gromov-Witten invariants where all $k_i$ are zero are called
\emph{primary}; they have an interpretation as the ``number'' (in a
suitable sense) of algebraic curves of genus $g$ and degree $\beta$ in $V$
which meet $n$ sufficiently generic cycles representing the Poincar\'e
duals of the cohomology classes $x_i$ (see Ruan \cite{Ruan}).

Gromov-Witten invariants in which some (or all) of the numbers $k_i$ are
positive are called \emph{descendent}; these do not admit so easily of an
enumerative interpretation. In genus $0$ and $1$, decsendent Gromov-Witten
invariants may be expressed in terms of the primary Gromov-Witten
invariants, by means of the topological recursion relations
\cite{KM1,genus2}. In higher genus, this is no longer the case: using
topological recursion relations, one can express genus $2$ Gromov-Witten
invariants in terms of those with $k_1+\dots+k_n\le1$ (see \cite{genus2}),
but no better. (It is likely that one can express genus $g$ descendent
Gromov-Witten invariants in terms of those with $k_1+\dots+k_n<g$.)

The most complete calculations of Gromov-Witten invariants have been made
for projective spaces. There are recursion relations among the primary
Gromov-Witten invariants in genus $0$ and in genus $1$ which determine
these invariants completely, and which follow respectively from the WDVV
equation and its analogue in genus $1$ (see \cite{elliptic}). But already
in genus $2$, the only known recursion relations involve both the primary
Gromov-Witten invariants and the descendent Gromov-Witten invariants with
$k_1+\dots+k_n=1$ (Belorusski-Pandharipande \cite{BP}).

Bearing the above facts in mind, it is not surprising that any conjecture
involving Gromov-Witten invariants in all genera, such as the Virasoro
conjecture, involves the consideration of descendent Gromov-Witten
invariants. We now turn to the formulation of this conjecture.

\subsection{The Novikov ring}
We employ the notation of Witten's foundational paper \cite{Witten}, except
that we explicitly introduce the Novikov ring $\Nov$ of $V$. Let
$H^+_2(V,\Z)$ denote the semi-group of $H_2(V,\Z)$ which is the image under
the cycle map of the semigroup of effective algebraic 1-cycles $\ZE_1(V)$
on $V$. The Novikov ring is
$$
\Nov = \biggl\{ a = \sum_{\beta\in H_2(V,\Z)} a_\beta q^\beta \Bigm| \text{
$\supp(a) \subset \beta_0 + H_2^+(V,\Z)$ for some $\beta_0\in H_2(V,\Z)$}
\biggr\} ,
$$
with product $q^{\beta_1}q^{\beta_2}=q^{\beta_1+\beta_2}$ and grading
$\deg\bigl(q^\beta\bigr)=-2c_1(V)\cap\beta$. The product on $\Lambda$ is
well-defined, since for any smooth projective variety $V$ with K\"ahler
form $\om$, the set $\{\beta\in H_2^+(V,\Z) \mid \om\cap\beta\le c\}$ is
finite for each $c>0$. By working over the Novikov ring, we may combine the
Gromov-Witten invariants in different degrees into a single generating
function:
$$
\<\tau_{k_1}(x_1)\dots\tau_{k_n}(x_n)\>^V_g = \sum_{\beta\in H_2^+(V,\Z)}
q^\beta \<\tau_{k_1}(x_1)\dots\tau_{k_n}(x_n)\>^V_{g,\beta} .
$$

\subsection{The large phase space $\PHASE(V)$}
If $V$ is a smooth projective variety, let $\{\gamma_a\mid a\in A\}$ be a
basis for $H(V)=H^\bull(V,\C)$; denote by $0\in A$ a distinguished index
with $\gamma_0=1\in H^0(V,\C)$. We suppose that this basis is homogeneous
with respect to the Hodge decomposition: each $\gamma_a$ is in
$H^{p_a,q_a}(V)$ for some $p_a$ and $q_a$.

Let $\phase(V)$ be the formal superscheme over $\Lambda$ obtained by
completing the affine superspace $H(V)$ at $0$; in the physics literature,
it is called the \emph{small phase space}. This formal superscheme has
coordinates $\{ u^a \mid a \in A \}$; denote the vector field $\p/\p u^a$
on $\phase(V)$ by $\p_a$. This is the superscheme on which the
Gromov-Witten invariants in genus $0$ define the structure of a Frobenius
supermanifold (Section \ref{frobenius}).

More important for us will be an infinite-dimensional formal superscheme
$\PHASE(V)$ defined over $\Lambda$ which is obtained by completing the
affine superspace
$$
H^\bull_{S^1}(V,\C) \cong H(V)[\om] , \quad \om \in H^2_{S^1} ,
$$
at $0$; physicists call this the \emph{large phase space}. This formal
superscheme has coordinates $\{ t_m^a \mid a\in A, m\ge0 \}$; denote the
vector field $\p/\p t^a_m$ on $\PHASE(V)$ by $\p_{m,a}$.

Note that if $H^{\text{odd}}(V,\C)=0$, then both $\phase(V)$ and
$\PHASE(V)$ are formal \emph{schemes} over $\Lambda$, since in that case
all coordinates have even $\Z/2$-grading.

In writing formulas in the coordinate systems $\{u^a\}$ and $\{t_m^a\}$, we
always assume the summation convention over indices $a$, $b, \dots\in A$,
using the non-degenerate inner product $\eta_{ab} = \int_V \gamma_a \cup
\gamma_b$ and its inverse $\eta^{ab}$ to raise and lower indices as needed.

\subsection{The Gromov-Witten potential}
The genus $g$ Gromov-Witten potential of $V$ is the function on the
superscheme $\PHASE(V)$ defined by the formula
$$
\<\<~\>\>^V_g = \sum_{n=0}^\infty \frac{1}{n!} \sum_{\substack{k_1\dots k_n
\\ a_1\dots a_n}} t_{k_n}^{a_n} \dots t_{k_1}^{a_1} \< \tau_{k_1,a_1} \dots
\tau_{k_n,a_n} \>^V_g ,
$$
where $\tau_{k,a}$ is an abbreviation for $\tau_k(\gamma_a)$. (The peculiar
ordering of the variables $t^{a_i}_{k_i}$ reflects the potential presence
of odd-dimensional cohomology classes on $V$.)

The total Gromov-Witten potential is
\begin{equation} \label{GW}
Z(V) = \exp \biggl( \sum_{g\ge 0} \hbar^{g-1} \<\< ~ \>\>_g^V \biggr) .
\end{equation}
This potential does not lie in any space of functions on $\PHASE(V)$:
rather, it defines a line bundle on $\PHASE(V)$, whose sections are objects
of the form
$$
\sum_{k=-\infty}^\infty \hbar^k f_k \* Z(V) , \quad f_k \in \CO_{\PHASE(V)}
.
$$
This line bundle has a flat connection, given by the formula
$$
\p_{m,a} \Bigl( \sum_{k=-\infty}^\infty \hbar^k f_k \* Z(V) \Bigr) = \sum_k
\hbar^k \Bigl( \p_{m,a} f_k + \sum_{g=0}^\infty \<\<\tau_{m,a}\>\>_g^V
f_{k-g+1} \Bigr) \* Z(V) .
$$
We will refrain from mentioning this line bundle again, and pretend that
$Z(V)$ is actually a function on $\PHASE(V)$.

\subsection{The statement of the conjecture}
In \cite{EHX}, Eguchi~et~al.\ introduce a sequence of differential
operators $L_k$, $k\ge-1$, on the formal superscheme $\PHASE(V)$ (or
rather, on the line bundle associated to the section $Z(V)$). To state the
formulas for the operators $L_k$, we need some more notation.

Let $R^b_a$ be the matrix associated to multiplication on $H(V)$ by the
first Chern class $c_1(V)$ (or equivalently, anticanonical class $-K_V$) of
$V$, defined by
$$
R^b_a \gamma_b = c_1(V) \cup \gamma_a .
$$
Let $\mu$ be the diagonal matrix with entries
$$
\mu_a = p_a - \tfrac{r}{2} ,
$$
where $r=\dim_\C(V)$. Let $[x]^k_i = e_{k+1-i}(x,x+1,\dots,x+k)$, where
$e_k$ is the $k$th elementary symmetric function of its arguments; thus,
$$
\sum_{i=0}^{k+1} s^i [x]^k_i = (s+x)(s+x+1)\dots(s+x+k) .
$$

The differential operators $L_k$ are defined by the following formula:
\begin{align} \label{Lk}
L_k = \sum_{i=0}^{k+1} & \biggl( \tfrac{\hbar}{2} \sum_{m=i-k}^{-1} (-1)^m
[\mu_a\!+\!m\!+\!\half]^k_i (R^i)^{ab} \p_{-m-1,a} \p_{m+k-i,b} \\ & -
[\tfrac{3-r}{2}]^k_i (R^i)^b_0 \p_{k-i+1,b} + \sum_{m=0}^\infty
[\mu_a\!+\!m\!+\!\half]^k_i (R^i)^b_a t^a_m \p_{m+k-i,b} \biggr) \notag \\
+ & \tfrac{1}{2\hbar} (R^{k+1})_{ab} t^a_0 t^b_0 + \tfrac{\delta_{k,0}}{48}
\int_V \bigl( (3-r) c_r(V) - 2c_1(V) c_{r-1}(V) \bigr) . \notag
\end{align}
(It is understood that $\p_{a,m}$ vanishes if $m<0$.) All of these
operators are quadratic expressions in the operators of multiplication by
$t_m^a$ and differentiation $\p_{m,a}$.

The formula for the first of these operators $L_{-1}$ is far simpler than
the others, and does not involve the coefficients $\mu_a$ and $R_a^b$:
\begin{equation} \label{L-1}
L_{-1} = - \p_{0,0} + \sum_{m=1}^\infty t^a_m \p_{m-1,a} + \frac{1}{2\hbar}
\eta_{ab} t^a_0 t^b_0 .
\end{equation}
This operator and $L_0$ are first-order differential operators, whereas
$L_k$ is a second-order differential operator for $k>0$.

In the original paper \cite{EHX}, the operators $L_k$ were only introduced
for $V$ a Grassmannian; the above extension may be found in \cite{EJX} and
is due to S. Katz. We may now state the Virasoro conjecture \cite{EHX,EJX}.
\begin{conjecture}
If $V$ is a smooth projective variety over $\C$,
$$
\text{$L_kZ(V)=0$ for all $k\ge-1$.}
$$
\end{conjecture}

The reason that this conjecture is called the Virasoro conjecture is that
the operators $L_k$ satisfy the commutation relations
\begin{equation} \label{virasoro}
[L_k,L_\ell] = (k-\ell) L_{k+\ell} ,
\end{equation}
and thus form a Lie subalgebra of the Virasoro algebra isomorphic to the
Lie algebra of polynomial vector fields on the line, with basis
$$
L_k = - \zeta^{k+1} \frac{\p}{\p\zeta} .
$$
In particular, the proof of the commutation relation $[L_1,L_{-1}]=2L_0$
depends on the Hirzebruch-Riemann-Roch theorem.

Granted the commutation relations \eqref{virasoro}, we have
$$
L_k = \frac{(-1)^{k-2}}{(k-2)!} \ad(L_1)^{k-2} L_2 , \quad k\ge2 ;
$$
together with the formula $L_1=-\frac13[L_{-1},L_2]$, this shows that the
Virasoro conjecture is equivalent to the formulas $L_{-1}Z(V)=0$ and
$L_2Z(V)=0$. However, this observation appears to be of little practical
importance in understanding the conjecture.

\subsection{What is known}
Very little headway has been made in the proof of the Virasoro
conjecture. In this section, we summarize what is presently known; we
present more details of all but the work of Dubrovin and Zhang \cite{DZ} in
genus $1$ in later sections.

The Virasoro conjecture in genus $0$ is a formal consequence of simple
properties of the Gromov-Witten invariants in genus $0$; we give a new
proof based on the sketch in \cite{EHX} in Section~\ref{genus0}.

The original Virasoro conjecture, in the special case where $V$ is a point,
was discovered by Dijkgraaf, Verlinde and Verlinde \cite{DVV}; their work
was one of the main influences which led to the formulation of the general
conjecture. They showed that in this case, the conjecture is equivalent to
Witten's conjecture \cite{Witten} relating the Gromov-Witten potential of a
point to the KdV hierarchy; a way to prove both of these conjectures was
given by Kontsevich \cite{K}.

Recently, Dubrovin and Zhang \cite{DZ} have conjectured that for any $V$
whose small phase space $\phase(V)$ is a semisimple Frobenius manifold (a
condition on the genus $0$ Gromov-Witten invariants of $V$, satisfied, for
example, if $V$ is a Grassmannian), there is a \emph{unique} hierarchy
compatible with the Virasoro conjecture in the same way that Witten's
conjecture is compatible with the Virasoro conjecture for a point. They
have constructed this hierarchy up to genus $1$, and used it to give a
proof of the Virasoro conjecture in genus $1$ for such varieties.

The other piece of evidence which led to the formulation of the Virasoro
conjecture is that the equations $L_{-1}Z(V)=0$ and $L_0Z(V)=0$ hold for
arbitrary $V$. The proofs, due respectively to Witten \cite{Witten} and
Hori \cite{Hori}, will be recalled in Sections \ref{Puncture} and
\ref{Hori}.

The Virasoro conjecture simplifies greatly for Calabi-Yau manifolds
(projective varieties for which $c_1(V)=0$ and $H^1(V,\C)$), since in that
case, the matrix $R^b_a$ vanishes, and, for $k>0$,
\begin{align*}
L_k &= -
\tfrac{\Gamma\bigl(\mu_a+\frac{5-r}{2}\bigr)}{\Gamma(\frac{3-r}{2})}
\p_{k+1,0} + \sum_{m=0}^\infty
\tfrac{\Gamma(\mu_a+m+k+\frac32)}{\Gamma(\mu_a+m+\frac12)} t^a_m \p_{m+k,a}
\\ & \qquad + \tfrac{\hbar}{2} \sum_{m=-k}^{-1} (-1)^m
\tfrac{\Gamma(\mu_a+m+k+\frac32)}{\Gamma(\mu_a+m+\frac12)} \eta^{ab}
\p_{-m-1,a} \p_{m+k,b} .
\end{align*}
We prove in Section \ref{CY} that the Virasoro conjecture holds in genus
$g>0$ for Calabi-Yau manifolds of dimension at least $3$, by purely
dimensional arguments. It follows that the conjecture yields no constraints
on the Gromov-Witten invariants in such dimensions.

We close with one last piece of ``evidence'' for the Virasoro
conjecture. For a general smooth projective variety $V$, we may extract the
coefficient of $q^0$ from the formula for $z_{k,g}$, and it turns out
\cite{GP} that the resulting equation depends on $V$ only through its
dimension $r$. The resulting equations are empty if $r>2$, but for curves
and surfaces, we obtain the following interesting implications of the
Virasoro conjecture (independent of the curve, respectively surface, in
question):

(1) for curves, if $2g+n-3=k_1+\dots+k_n$, then
$$
\int_{\Mbar_{g,n}} \Psi_1^{k_1} \dots \Psi_n^{k_n} \lambda_g =
\frac{(2g+n-3)!}{k_1! \dots k_n!} \int_{\Mbar_{g,1}} \Psi_1^{2g-2}
\lambda_g .
$$

(2) for surfaces, $g+n-2=k_1+\dots+k_n$ and $k_i>0$,
$$
\int_{\Mbar_{g,n}} \Psi_1^{k_1} \dots \Psi_n^{k_n} \lambda_g \lambda_{g-1}
= \frac{(2g-1)!! (2g+n-3)!}{(2g-1)! (2k_1-1)!! \dots (2k_n-1)!!}
\int_{\Mbar_{g,1}} \Psi_1^{g-1} \lambda_g \lambda_{g-1} .
$$

Remarkably, the second of these formulas had earlier been conjectured, in
an entirely different context, by Faber \cite{Faber}; he has proved it in
genera up to $15$.

\subsection*{Acknowledgements}
I am grateful to the organizers of ``Hirzebruch 70'' and the Banach
Institute and to J.-P. Bismut of Universit\'e de Paris-Sud for inviting me
to lecture on the Virasoro conjecture, and to E. Arbarello and la Scuola
Normale Superiore, Pisa, where a large part of this review was written.

In learning this subject, I profitted greatly from conversations with
J. Bryan, T. Eguchi, E. Frenkel, S. Katz, C.-S. Xiong, and from
collaboration with R. Pandharipande.

\section{The definition of Gromov-Witten invariants} \label{definition}

\subsection{Stable maps}
The Gromov-Witten invariants of a projective manifold reflect the
intersection theory of the moduli spaces of stable maps
$\Mbar_{g,n}(V,\beta)$, whose definition we now recall.

Let $V$ be a smooth projective variety. (In this paper, all varieties are
defined over $\C$.) A \emph{prestable map}
$$\textstyle
(f:C\to V,z_1,\dots,z_n)
$$
of genus $g\ge0$ and degree $\beta\in H_2^+(V,\Z)$ with $n$ marked points
consists of the following data:
\begin{enumerate}
\item a connected projective curve $C$ of arithmetic genus
$g=h^1(C,\CO_C)$, whose only singularities are ordinary double points,
\item $n$ distinct smooth points $(z_1,\dots,z_n)$ of $C$;
\item an algebraic map $f:C\to V$, such that the cycle $f_*[C]\in
H^2(V,\Z)$ equals $\beta$.
\end{enumerate}
If $\tilde{C}$ is the normalization of $C$, the \emph{special points} in
$\tilde{C}$ are the inverse images of the singular and marked points of
$C$. (Note that the degree of $f:C\to V$ equals $0$ if and only if its
image is a single point.)

A prestable map $(f:C\to V,z_1,\dots,z_n)$ is \emph{stable} if it has no
infinitesimal automorphisms fixing the marked points. The condition of
stability is equivalent to the following: each irreducible component of
$\tilde{C}$ of genus $0$ on which $f$ has degree $0$ has at least $3$
special points, while each irreducible component of $\tilde{C}$ of genus
$1$ on which $f$ has degree $0$ has at least $1$ special point. In
particular, there are no stable maps of genus $g$ and degree $0$ with $n$
marked points unless $2(g-1)+n>0$.

\subsection{The moduli stack $\Mbar_{g,n}(V,\beta)$ of stable maps}
Let $\Mbar_{g,n}(V,\beta)$ be the moduli stack of $n$-pointed stable maps
of genus $g$ and degree $\beta$, introduced by Kontsevich
\cite{KM}. Behrend and Manin \cite{BM} show that $\Mbar_{g,n}(V,\beta)$ is
a proper Deligne-Mumford stack (though not in general smooth). When $n=0$,
we write $\Mbar_g(V,\beta)$ instead of $\Mbar_{g,n}(V,\beta)$.

An important role in the theory is played by the map
\begin{equation} \label{universal}
\pi : \Mbar_{g,n+1}(V,\beta) \to \Mbar_{g,n}(V,\beta) .
\end{equation}
This is the operation which forgets the last point $z_{n+1}$ of a stable
map $(f:C\to V,z_1,\dots,z_{n+1})$, leaving a prestable map $(f:C\to
V,z_1,\dots,z_n)$, and contracts any rational component of $C$ on which $f$
has zero degree and which obstructs the stability of $(f:C\to
V,z_1,\dots,z_n)$. Behrend and Manin show that $\pi$ is a flat map, whose
fibre at $(f:C\to V,z_1,\dots,z_n)$ may be identified in a natural way with
the curve $C$; this identifies $\Mbar_{g,n+1}(V,\beta)$ with the universal
curve $\Cbar_{g,n}(V,\beta)$ over $\Mbar_{g,n}(V,\beta)$.

\subsection{The virtual fundamental class}
Denote by $f:\Cbar_{g,n}(V,\beta)\to V$ the universal stable map, defined
by sending the stable map $(f:C\to V,z_1,\dots,z_{n+1})$ to
$f(z_{n+1})$. If the sheaf $R^1\pi_*f^*TV$ vanishes, the
Grothendieck-Riemann-Roch theorem implies that the stack
$\Mbar_{g,n}(V,\beta)$ is smooth, of dimension
\begin{equation} \label{vdim}
\vdim\Mbar_{g,n}(V,\beta) = (3-r)(g-1) + c_1(V)\cap\beta + n .
\end{equation}
In general, we call this number the \emph{virtual dimension} of
$\Mbar_{g,n}(V,\beta)$.

The hypothesis $R^1\pi_*f^*TV=0$ is rarely true. However, there is an
algebraic cycle
$$
[\Mbar_{g,n}(V,\beta)]^\virt \in
H_{2\vdim\Mbar_{g,n}(V,\beta)}(\Mbar_{g,n}(V,\beta),\Q) ,
$$
the \emph{virtual fundamental class}, which stands in for
$[\Mbar_{g,n}(V,\beta)]$ in the general case; this cycle is constructed in
Behrend \cite{B}, and in Li and Tian \cite{LT}. It is the existence of this
cycle which gives rise to the Gromov-Witten invariants.

One of the main properties of the virtual fundamental class is the
following formula (Axiom IV, Behrend \cite{B}):
\begin{equation} \label{flat}
\pi^*[\Mbar_{g,n}(V,\beta)]^\virt = [\Mbar_{g,n+1}(V,\beta)]^\virt .
\end{equation}
In particular, if $\vdim\Mbar_g(V,\beta)<0$, we see that
$[\Mbar_{g,n}(V,\beta)]^\virt=0$ for all $n\ge0$.

\subsection{Gromov-Witten invariants}
The projection \eqref{universal} has $n$ canonical sections
$$
\sigma_i : \Mbar_{g,n}(V,\beta)\to\Cbar_{g,n}(V,\beta) ,
$$
corresponding to the $n$ marked points of the curve $C$. Let
$$
\om=\om_{\Cbar_{g,n}(V,\beta)/\Mbar_{g,n}(V,\beta)}
$$
be the relative dualizing sheaf; the line bundle $\Om_i = \sigma_i^*\om$
has fibre $T^*_{z_i}C$ at the stable map $(f:C\to V,z_1,\dots,z_n)$. Let
$\Psi_i$ be the cohomology class $c_1(\Om_i)$.

Let $\ev:\Mbar_{g,n}(V,\beta)\to V^n$ be evaluation at the marked points:
$$
\ev : { \textstyle (f:C\to V,z_1,\dots,z_n) } \mapsto \bigl(
f(z_1),\dots,f(z_n) \bigr) \in V^n .
$$
If $x_1,\dots,x_n$ are cohomology classes of $V$, we define the
Gromov-Witten invariants by the formula
$$
\< \tau_{k_1}(x_1) \dots \tau_{k_n}(x_n) \>^V_{g,\beta} =
\int_{[\Mbar_{g,n}(V,\beta)]^\virt} \Psi_1^{k_1} \dots \Psi_n^{k_n} \cup
\ev^*(x_1\boxtimes \dots \boxtimes x_n) .
$$
In addition to the generating functions $\<\<~\>\>_g^V$ considered in the
introduction, we will also work with the functions
$$
\<\<\tau_{k_1,a_1}\dots\tau_{k_n,a_n}\>\>^V_g =
\p_{k_1,a_1}\dots\p_{k_n,a_n} \<\<~\>\>_g^V \in \CO_{\PHASE(V)} ,
$$
which equal $\<\tau_{k_1,a_1}\dots\tau_{k_n,a_n}\>^V_g$ at $0\in
\PHASE(V)$.

\subsection{Puncture equation} \label{Puncture}
Witten \cite{Witten} proves the following equations:
\begin{equation} \label{puncture}
\<\tau_{0,0}\tau_{k_1,a_1}\dots\tau_{k_n,a_n}\>_{g,\beta}^V = \sum_{i=1}^n
\<\tau_{k_1,a_1}\dots\tau_{k_i-1,a_i} \dots \tau_{k_n,a_n}\>_{g,\beta}^V .
\end{equation}
In degree zero, there is one exceptional case:
\begin{equation} \label{eta}
\<\tau_{0,0}\tau_{0,a}\tau_{0,b}\>_{0,0}^V = \eta_{ab} .
\end{equation}
These equations are a simple consequence of the geometry of the divisors
associated to the line bundles $\Om_i$ (cf.\ \cite{KM1,genus2}), combined
with \eqref{flat}.

Together, \eqref{puncture} and \eqref{eta} are equivalent to the sequence
of equations
$$
\<\<\tau_{0,0}\>\>_g^V = \sum_{m=1}^\infty t_m^a \<\<\tau_{m-1,a}\>\>_g^V +
\tfrac{1}{2} \delta_{g,0} \eta_{ab} t_0^a t_0^b .
$$
We may combine these into a single equation by multiplying by $\hbar^{g-1}$
and summing over $g$:
$$
\biggl( - \p_{0,0} + \sum_{m=1}^\infty t_m^a \p_{m-1,a} \biggr) \biggl(
\sum_{g=0}^\infty \hbar^{g-1} \<\<~\>\>_g^V \biggr) + \frac{1}{2\hbar}
\eta_{ab} t_0^a t_0^b = 0 .
$$
In terms of the Gromov-Witten potential $Z(V)$, this becomes a homogeneous
first-order linear differential equation, known as the \emph{puncture
equation} (or alternatively, the string equation):
$$
\biggl( - \p_{0,0} + \sum_{m=1}^\infty t_m^a \p_{m-1,a} + \frac{1}{2\hbar}
\eta_{ab} t_0^a t_0^b \biggr) Z(V) = 0 .
$$
The differential operator on the left-hand side of this equation is
precisely the operator $L_{-1}$ of \eqref{L-1}; thus, the puncture equation
$L_{-1}Z(V)=0$ is actually a part of the Virasoro conjecture.

\subsection{Divisor equation}
If $\om\in H^2(V,\C)$, let $R_a^b(\om)$ be the matrix of multiplication by
$\om$ on $H(V)$: $\om \cup \gamma_a = R_a^b(\om) \gamma_b$. By the same
method as the puncture equation is proved, Hori \cite{Hori} proves the
\emph{divisor equation}
\begin{align} \label{divisor}
\<\tau_0(\om)\tau_{k_1,a_1}\dots\tau_{k_n,a_n}\>_{g,\beta}^V &= \bigl( \om
\cap \beta \bigr) \* \<\tau_{k_1,a_1}\dots\tau_{k_n,a_n}\>_{g,\beta}^V \\
\notag & + \sum_{i=1}^n R_{a_i}^b(\om)
\<\tau_{k_1,a_1}\dots\tau_{k_i-1,b} \dots \tau_{k_n,a_n}\>_{g,\beta}^V .
\end{align}
In degree zero, there are two exceptional cases:
$$
\<\tau_0(\om)\tau_{0,a}\tau_{0,b}\>_{0,0}^V = R_{ab}(\om)
\quad\text{and}\quad \<\tau_0(\om)\>_{1,0}^V = \frac{1}{24} \int_V \om \cup
c_{r-1}(V) .
$$

\subsection{Dilaton equation}
Witten also proves the following equations in \cite{Witten}:
$$
\<\tau_{1,0}\tau_{k_1,a_1}\dots\tau_{k_n,a_n}\>_g^V = (2g-2+n)
\<\tau_{k_1,\alpha_1}\dots\tau_{k_n,\alpha_n}\>_g^V .
$$
In degree zero, there is one exceptional case:
$$
\<\tau_{1,0}\>_{1,0}^V = \frac{\chi(V)}{24} .
$$
As in the discussion of the puncture equation, these equations may be
combined into a first-order differential equation
$$
\left( \Dil + \frac{\chi(V)}{24} \right) Z(V) = 0 ,
$$
called the \emph{dilaton equation}, where $\Dil$ is the differential
operator
\begin{equation} \label{Dilaton}
\Dil = - \p_{1,0} + \sum_{m=0}^\infty t_m^a \p_{m,a} + 2 \hbar
\frac{\p}{\p\hbar} .
\end{equation}

\subsection{A characteristic number}
Let us introduce an abbreviation for the constant term of $L_0$ in
\eqref{Lk}:
$$
\rho(V) = \frac{1}{48} \int_V \bigl( (3-r) c_r(V) - 2c_1(V) c_{r-1}(V)
\bigr) .
$$
This characteristic number behaves as follows under products:
$$
\rho(V\times W) = \rho(V)\chi(W) + \chi(V)\rho(W) - \tfrac{1}{16}
\chi(V)\chi(W) .
$$
Note that $\rho(V)$ vanishes for Calabi-Yau threefolds, hinting at the
special role which they play in the theory.

\subsection{Hori's equation} \label{Hori}
Because the dilaton operator $\Dil$ involves differentiation with respect
to the parameter $\hbar$, it is not a vector field on the large phase space
$\PHASE(V)$. By judiciously combining it with the equations for the virtual
dimension of the moduli spaces $\Mbar_{g,n}(V,\beta)$ and with the divisor
equation, Hori \cite{Hori} was able to construct from it a first-order
differential operator on $\PHASE(V)$ which annihilates the Gromov-Witten
potential $Z(V)$. As in the introduction, let us write $R_a^b$ for the
matrix $R_a^b(c_1(V))$.
\begin{theorem}
We have $L_0Z(V)=0$, where $L_0$ is the differential operator
\begin{align*}
L_0 &= - \half (3-r) \p_{1,0} +\sum_{m=0}^\infty (\mu_a\!+\!m\!+\!\half)
t_m^a \p_{m,a} - R_0^b \p_{0,b} + \sum_{m=1}^\infty R_a^b t_m^a \p_{m-1,b}
\\ & \quad + \frac{1}{2\hbar} R_{ab} t^a_0 t^b_0 + \delta_{k,0} \rho(V) .
\end{align*}
\end{theorem}
\begin{proof}
The formula \eqref{vdim} for the dimension of the virtual fundamental class
of $\Mbar_{g,n}(V,\beta)$ implies the following equations among
Gromov-Witten invariants:
$$
\sum_{i=1}^n (p_{\alpha_i}\!+\!k_i-1)
\<\tau_{k_1,a_1}\dots\tau_{k_n,a_n}\>_{g,\beta}^V =
\vdim\Mbar_g(V,\beta) \*
\<\tau_{k_1,a_1}\dots\tau_{k_n,a_n}\>_{g,\beta}^V
$$
In order to eliminate the dependence on the degree $\beta$, we subtract the
divisor equation \eqref{divisor} with $\om=c_1(V)$, obtaining the
differential equation
\begin{multline*}
\biggl( \sum_{m=0}^\infty (p_a\!+\!m-1) t^a_m \p{}_{m,a} - R_0^b \p{}_{0,b}
+ \sum_{m=1}^\infty R_a^b t^a_m \p{}_{m-1,a} \\ + \frac{1}{2\hbar} R_{ab}
t^a_0 t^b_0 - \frac{1}{24} \int_V c_1(V) \cup c_{r-1}(V) + (r-3) \hbar
\frac{\p}{\p\hbar} \biggr) Z(V) = 0 .
\end{multline*}
It only remains to eliminate the coefficient of $\p/\p\hbar$ in this
operator, by adding $\half(3-r)$ times the dilaton equation.
\end{proof}

The operators $L_{-1}$ and $L_0$ satisfy the commutation relation
$[L_0,L_{-1}]=L_{-1}$. Motivated by this relation, Eguchi~et~al.\ were led
to introduce the sequence of differential operators $L_k$, $k>0$, of
\eqref{Lk}.
\begin{theorem} \label{Virasoro}
The sequence of differential operators $L_k$, $k>0$, of \eqref{Lk} satisfy
the Virasoro commutation relations \eqref{virasoro}.
\end{theorem}

Before giving the proof, we need a small amount of quantum field
theory. Let $\Psi$ be the space of pseudodifferential operators on the
affine line with coordinate $z$, that is, expressions of the form
$$
\PP = \sum_{n=-\infty}^\infty p_n(z)\p^n ,
$$
where $p_n(z)\in\C[z]$. Let $\phi(z)$ be the generating
function (or \emph{free field})
\begin{equation} \label{free}
\phi^a(z) = \sum_{m=0}^\infty \frac{\Gamma(\half)}{\Gamma(m+\frac32)}
z^{m+\frac12} t_m^a - \hbar \sum_{m=0}^\infty
\frac{\Gamma(m+\frac12)}{\Gamma(\half)} z^{-m-\frac12} \eta^{ab}
\p_{m,b} - \tfrac{4}{3} \delta_0^a z^{3/2} .
\end{equation}
\begin{definition}
A differential operator $\delta$ on the large phase space $\PHASE(V)$ has
\emph{symbol} $\sigma(\delta)\in\Psi\o\End(H(V))$ if
$\sigma(\delta)\phi(z) + [\delta,\phi(z)] = 0$.
\end{definition}

This definition is justified by the following two lemmas.
\begin{lemma}
If $\sigma(\delta)=0$, then $\delta$ is a multiple of the identity
operator.
\end{lemma}
\begin{proof}
If $\sigma(\delta)=0$, then $\delta$ must commute with all of the
coefficients of the fields $\phi^a(z)$. Any such operator lies in the
center of the algebra of differential operators, whence the lemma.
\end{proof}

\begin{lemma}
The symbol of $[\sigma(\delta_1),\sigma(\delta_2)]$ equals
$[\delta_1,\delta_2]$.
\end{lemma}
\begin{proof}
$[\sigma(\delta_1),\sigma(\delta_2)] \phi(z) +
[[\delta_1,\delta_2],\phi(z)]$
\begin{align*}
&= \sigma(\delta_1) \sigma(\delta_2) \phi(z) - \sigma(\delta_2)
\sigma(\delta_1) \phi(z) + [\delta_1,[\delta_2,\phi(z)]] -
[\delta_2,[\delta_1,\phi(z)]] \\
& = - \sigma(\delta_1) [\delta_2,\phi(z)] + \sigma(\delta_2)
[\delta_1,\phi(z)] + [\delta_1,[\delta_2,\phi(z)]] -
[\delta_2,[\delta_1,\phi(z)]] \\
&= - [\delta_2,\sigma(\delta_1)\phi(z)+[\delta_1,\phi(z)]] + 
[\delta_1,\sigma(\delta_2)\phi(z)+[\delta_2,\phi(z)]] = 0 .
\qed\end{align*}
\def\qed{}
\end{proof}

The reader may have observed that the formulas for the field $\phi(z)$ and
the operators $L_k$ and $\Dil$ simplify if we rewrite them in terms of the
shifted coordinates
$$
\tilde{t}_m^a = \begin{cases} t_1^0 - 1 , & m=1, a=0, \\
t_m^a , & \text{otherwise.}
\end{cases}$$
For example, the formulas for $L_{-1}$ and $\Dil$ become simply
$$
L_{-1} = \sum_{m=1}^\infty \tilde{t}^a_m \p_{m-1,a} + \frac{1}{2\hbar}
\eta_{ab} \tilde{t}^a_0 \tilde{t}^b_0 \quad,\quad \Dil = \sum_{m=0}^\infty
\tilde{t}_m^a \p_{m,a} + 2 \hbar \frac{\p}{\p\hbar} ,
$$
while that for $\phi(z)$ becomes
$$
\phi^a(z) = \sum_{m=0}^\infty \frac{\Gamma(\half)}{\Gamma(m+\frac32)}
z^{m+\frac12} \tilde{t}_m^a - \hbar \sum_{m=0}^\infty
\frac{\Gamma(m+\frac12)}{\Gamma(\half)} z^{-m-\frac12} \eta^{ab}
\p_{m,b} .
$$
The explanation of this is as follows: consider the Lie superalgebra $\g$
of differential operators on $\PHASE(V)$ quadratic in the operators
$\p_{m,a}$ and $t_m^a$. This Lie superalgebra has an increasing sequence of
subspaces, define inductively by
$$
F_k\g = \{ A\in\g \mid [L_{-1},A]\in F_{k-1}\g \} \text{, where}\quad
F_{-1}\g=\<1,L_{-1}\> .
$$
By induction on $k$, we may show that $[F_k\g,F_\ell\g]\subset
F_{k+\ell}\g$, and hence that the union
$$
F_\infty\g = \bigcup_{k=0}^\infty F_k\g
$$
is a Lie sub-superalgebra of $\g$. Let $\DD=z\p$. One may show
\cite{psi} that every element of $F_\infty\g$ has a symbol of the form
\begin{equation} \label{f(D)}
\sigma(\delta) = \sum_{k=0}^\infty f_k(\DD) \p^k \in \Psi\o\End(H(V)) ,
\end{equation}
where $f_k(\DD)\in\End(H(V))[\DD]$ satisfies
$f_k(-t)=(-1)^{k+1}f_k^*(t-k)$. The space of such pseudodifferential
operators forms a Lie superalgebra, of which $F_\infty\g$ is a central
extension, associated to the two-cocycle
\begin{multline} \label{cocycle}
c(f(\DD)\p^k,g(\DD)\p^\ell) = \frac{\delta_{k+\ell,0}}{2} \\ \Bigl(
\sum_{m=0}^{k-1} \Str\bigl( f(m-k+\tfrac{1}{2}) g(m+\tfrac12) \bigr) -
\sum_{m=0}^{-k-1} \Str\bigl( f(m+\tfrac{1}{2}) g(m-k+\tfrac12) \bigr)
\Bigr) ,
\end{multline}
where $\Str$ is the supertrace (the difference of the traces on the even
and odd degree subspaces of $H(V)$). This two-cocyle is calculated using
the natural section $\sigma^{-1}$ of the symbol map
\begin{align} \label{sigma-1}
\sigma^{-1}\bigl( f(\DD) \p^k \bigr) &= \frac{\hbar}{2} \sum_{m=0}^{-k-1}
(-1)^m \eta^{ab} f\bigl(m+\tfrac12\bigr){}_a^b \p_{m,c} \p_{-k-m-1,b} \\
\notag & - \sum_{m=(-k){}_+}^\infty f\bigl(m+\tfrac{1}{2}\bigr){}_a^b
\tilde{t}_{m+k}^a \p_{m,b} \\ \notag & + \frac{1}{2\hbar} \sum_{m=0}^{k-1}
(-1)^{m+1} \eta_{bc} f\bigl(m-k+\tfrac{1}{2}\bigr){}_a^b
\tilde{t}_m^a\tilde{t}_{k-m-1}^c .
\end{align}

\begin{proof}[Proof of Theorem \ref{Virasoro}]
Let $\mu$ be the diagonal matrix $\mu_a^b=\delta_a^b\mu_a$, and let $R$ be
the matrix with entries $R_a^b$. Then
\begin{equation} \label{lk}
L_k = - \sigma^{-1} \Bigl( (z+\mu\,\p^{-1}+R)^{k+1} \p \Bigr) +
\delta_{k,0} \rho(V) \in F_k\g\subset F_\infty\g .
\end{equation}
Since $\sigma(L_k)= -(z+\mu\,\p^{-1}+R)^{k+1} \p$, we see that
$$
[\sigma(L_k),\sigma(L_\ell)] = (k-\ell) \sigma(L_{k+\ell}) ,
$$
and hence that
$$
[L_k,L_\ell] = (k-\ell) L_{k+\ell} + c(\sigma(L_k),\sigma(L_\ell)) - 2
\delta_{k+\ell,0} \rho(V) ,
$$
where $c$ is the two-cocyle of \eqref{cocycle}.

It remains to show that
\begin{equation} \label{k+l}
c(\sigma(L_k),\sigma(L_\ell)) = 2\delta_{k+\ell,0} \rho(V) .
\end{equation}
We have
$$
c(\sigma(L_k),\sigma(L_\ell)) = c\Bigl( (z+\mu\,\p^{-1}+R)^{k+1} \p ,
(z+\mu\,\p^{-1}+R)^{\ell+1} \p \Bigr) .
$$
Since the matrix $R$ raises degree while $z+\mu\,\p^{-1}$ preserves
degree, terms involving $R$ do not contribute to the supertrace, showing
that
$$
c(\sigma(L_k),\sigma(L_\ell)) = c\Bigl( (z+\mu\,\p^{-1})^{k+1} \p ,
(z+\mu\,\p^{-1})^{\ell+1} \p \Bigr) .
$$
Since $(z+\mu\,\p^{-1})^{k+1} \p = (\DD+\mu)(\DD+\mu-1)\dots(\DD+\mu-k)
\p^{-k}$, the formula \eqref{cocycle} for the cocycle $c$ shows that
$c(\sigma(L_k),\sigma(L_\ell))$ vanishes unless $k+\ell=0$ and $|k|=1$, and
that
$$
c(\sigma(L_1),\sigma(L_{-1})) = c((\DD+\mu)(\DD+\mu-1)\p^{-1} , \p ) = -
\tfrac{1}{2} \Str\bigl( \mu^2-\tfrac14 \bigr) .
$$
Eq.\ \eqref{k+l} is now a consequence of the following formula of Libgober
and Wood \cite{LW} (cf.\ Borisov \cite{Borisov}).
\end{proof}

\begin{proposition} \label{Libgober}
$$
\Str\bigl( \mu^2 \bigr) = \frac{1}{12} \int_V \bigl( rc_r(V) + 2 c_1(V)
c_{r-1}(V) \bigr)
$$
\end{proposition}
\begin{proof}
Let $\chi_t(V) = \sum_{p=0}^r t^p \chi(V,\Om^p)$ be the Hirzebruch
characteristic, and write $h(t)=\chi_{-t}(V)$. We have
\begin{align*}
\Str\bigl(\mu^2\bigr) &= \sum_{p=0}^r \bigl(p-\tfrac{r}{2}\bigr)^2 (-1)^p
\chi(V,\Om^p) = \left( t\frac{d}{dt}-\frac{r}{2} \right)^2 h(1) \\ &=
\ddot{h}(1) + (1-r) \dot{h}(1) + \frac{r^2}{4} h(1) .
\end{align*}
By the Hirzebruch-Riemann-Roch theorem,
$$
h(t) = \int_V \Todd(V) \sum_{p=0}^r (-t)^p \ch(\Lambda^pT^*V) .
$$
Introducing the Chern roots $c(V) = (1+x_1) \dots (1+x_r)$ of $V$, we see
that
\begin{align*}
h(1+t) &= \int_V \prod_{i=0}^r \frac{x_i}{1-e^{-x_i}} \sum_{p=0}^r (-1)^p
(1+t)^p \sum_{i_1<\dots<i_p} e^{-x_{i_1}-\dots-x_{i_p}} \\ &= \int_V
\prod_{i=0}^r \frac{x_i}{1-e^{-x_i}} ( 1 - (1+t)e^{-x_i} ) = \int_V
\prod_{i=0}^r \left( x_i - t \sum_{k=0}^\infty \frac{B_k}{k!} x_i^k \right)
.
\end{align*}
We conclude that
\begin{gather*}
h(1) = \int_V c_r(V) , \quad \dot{h}(1) = - \frac{r}{2} \int_V c_r(V) , \\
\ddot{h}(1) = \frac{r(r-1)}{4} \int_V c_r(V) + \frac{1}{6} \int_V \bigl(
c_1(V)c_{r-1}(V) - r c_r(V) \bigr) ,
\end{gather*}
and the result follows.
\end{proof}

\subsection{Remarks on the above proof}
Observe that the proof of Theorem \ref{Virasoro} had two parts:
\begin{enumerate}
\item showing that $[L_k,L_\ell]-(k-\ell)L_{k+\ell}$ is a constant --- this
only required of the matrices $\mu$ and $R$ that $[\mu,R]=R$ (that is, that
$\mu$ defines a grading of $H(V)$ in which $R$ raises degree by $1$);
\item showing that this constant vanishes --- this required Proposition
\ref{Libgober} to hold, a far more restrictive condition.
\end{enumerate}
There is no natural definition of a grading operator $\mu$ on compact
symplectic manifolds satisfying these conditions; this suggests,
although of course it does not prove, that there is no generalization
of the Virasoro conjecture to compact symplectic
manifolds\footnote{The equations $L_{-1}Z(M)=L_0Z(M)=0$ and the
restriction to genus $0$ of the conjecture do nevertheless hold for
any compact symplectic manifold $M$, provided that we take $\mu$ to
equal multiplication by $\frac14(2n-\dim_\R M)$ on $H^n(M,\C)$; the
obstruction is at genus $1$.}.

If we replace the holomorphic degree $p_a$ in the definition of the matrix
$\mu$ by the anti-holomorphic degree $q_a$, we obtain a new grading
operator $\bar{mu}$; the condition $[\mu,C]=C$ is satisfied by any affine
combination
\begin{equation} \label{mus}
\mu_s = (1-s) \mu + s \bar{\mu} ;
\end{equation}
let $L_k^s$ be the modified Virasoro operators obtained on replacing $\mu$
by $\mu^s$ in \eqref{Lk}. Following Borisov \cite{Borisov}, we have
$\Str\bigl( \mu_s^2 \bigr) = \Str\bigl( \mu^2 \bigr) +
s(s-1)\tilde{\rho}(V)$, where $\tilde\rho(V)=\Str\bigl( (\mu-\bar\mu)^2
\bigr)$, and hence
$$ \textstyle
[L^s_1,L_{-1}] = 2L^s_0 - \binom{s}{2} \tilde{\rho}(V) .
$$
If $\tilde\rho(V)\ne0$ and $s\ne0,1$, the modified Virasoro conjecture
$L_k^sZ(V)=0$ cannot hold.

Observe that $\tilde{\rho}(V)\ne0$ for smooth complete intersections of
sufficiently high degree (in particular, curves of nonzero genus), and for
Calabi-Yau threefolds. On the other hand, the formula $\tilde\rho(V\times
W) = \tilde\rho(V) \chi(W) + \chi(V) \tilde\rho(W)$ shows that
$\tilde{\rho}(V)=0$ for abelian varieties of dimension $r>1$.

It is curious that, if $r$ is even, there is an extension of the definition
of the operators $L_k$ to all $k\in\Z$, given by the formula \eqref{lk}
suitably interpreted. These operators satisfy the Virasoro relations with
central charge $\chi(V)$ (by \eqref{cocycle}).

\section{A basic lemma} \label{Key}

Let $z_{k,g}$ be the coefficient of $\hbar^{g-1}$ in $Z(V)^{-1}L_kZ(V)$.
The equation $L_kZ(V)=0$ is equivalent to the vanishing of $z_{k,g}$ for
all $g$. The explicit formula
\begin{align} \label{zkg}
z_{k,g} &= \sum_{i=0}^{k+1} \biggl( - [\tfrac{3-r}{2}]^k_i (R^i)^b_0
\<\<\tau_{k-i+1,b}\>\>_g^V + \sum_{m=0}^\infty [\mu_a\!+\!m\!+\!\half]^k_i
(R^i)^b_a t^a_m \<\<\tau_{m+k-i,b}\>\>_g^V \\ \notag &+ \frac{1}{2}
\sum_{m=i-k}^{-1} (-1)^m [\mu_a\!+\!m\!+\!\half]^k_i (R^i)^{ab} \Bigl(
\<\<\tau_{-m-1,a}\tau_{m+k-i,b}\>\>_{g-1}^V \\ \notag &+ \sum_{h=0}^g
\<\<\tau_{-m-1,a}\>\>_h^V \<\<\tau_{m+k-i,b}\>\>_{g-h}^V \Bigr) \biggr) +
\frac{\delta_{g,0}}{2} (R^{k+1})_{ab} t^a_0 t^b_0 +
\delta_{k,0}\delta_{g,1} \rho(V)
\end{align}
shows that $z_{k,g}$ depends on $\<\<~\>\>_{h}^V$ only for $h\le g$; in
particular, it is meaningful to speak of the Virasoro conjecture holding up
to genus $g$.

In this section, we show how the puncture equation, together with the
Virasoro relations \eqref{virasoro}, permits one to prove the Virasoro
conjecture in a given genus, provided we know that for all $k>0$, there is
an $i\ge0$ such that $\p_{0,0}^iz_{k,g}$ vanishes. This method lies behind
our proof of the Virasoro conjecture in genus $0$, and we expect it to be
equally useful in other situations.

Let $\CL_{-1}$ be the vector field part of the differential operator
$L_{-1}$:
$$
\CL_{-1}f = Z(V)^{-1} L_{-1} \bigl( fZ(V) \bigr) = - \p_{0,0}f +
\sum_{m=1}^\infty t_m^a \p_{m-1,a}f .
$$
This vector field has the following remarkable property. (The elegant
proof was provided by E. Frenkel.)
\begin{lemma} \label{key}
If $\p_{0,0}f$ and $\CL_{-1}f$ are constant, then so is $f$.
\end{lemma}
\begin{proof}
Let $E$ be the vector field
$$
E = \p_{0,0}+\CL_{-1} = \sum_{m=0}^\infty t_{m+1}^a \p_{m,a} .
$$
We must prove that if $Ef$ is constant, then so is $f$.

Together with $E$, the vector fields
$$
F = \sum_{m=0}^\infty (m+1)(m+2) t_m^a \p_{m+1,a}
$$
and
$$
H = \half [F,E] = \sum_{m=0}^\infty (m+1) t_m^a \p_{m,a} ,
$$
realize the Lie algebra $\Sl(2)$. It suffices to prove the lemma for
eigenfunctions of $H$, which are polynomial; on these, $F$ is locally
nilpotent, since the spectrum of $H$ is $\N$.

Suppose $f$ satisfies the hypotheses of the lemma, so that
$Ef=0$. Since $F^if=0$ for $i\gg0$, we see that the irreducible
$\Sl(2)$-module spanned by $f$ is finite-dimensional. But a
finite-dimensional representation of $\Sl(2)$ on which $H$ has
non-negative spectrum is a sum of trivial representations, and we
conclude that $Hf=0$, and hence that $f$ is constant.
\end{proof}

The above lemma is actually closely related to a result in formal
variational calculus: the algebra of polynomials $\Q[t_m^a\mid m\ge0,
a\in A]$ is isomorphic to the algebra of differential polynomials
$$
\Q\{u_a\mid a\in A\} = \Q[u_a^{(m)}\mid m\ge0, a\in A]
$$
under the identification of $t_m^a$ and $u_a^{(m)}$, and under this
isomorphism, the operator $\p_{0,0}+\CL_{-1}$ is carried into the
derivation
$$
\p = \sum_{m=0}^\infty u_a^{(m+1)} \frac{\p}{\p u_a^{(m)}} .
$$
Lemma \ref{key} shows that $\ker(\p)=\Q$, a result due to Gelfand and Dikii
(Section I.1, \cite{GD}).

\begin{lemma} \label{downward}
$\CL_{-1}z_{k,g}=-(k+1)z_{k-1,g}$
\end{lemma}
\begin{proof}
Since $[L_{-1},L_k]=-(k+1)L_{k-1}$ and $L_{-1}Z(V)=0$, we see that
$$
\CL_{-1}(Z(V)^{-1}L_kZ(V)) = Z(V)^{-1}L_{-1}L_kZ(V) = - (k+1) Z(V)^{-1}
L_{k-1} Z(V) .
$$
Extracting the coefficient of $\hbar^{g-1}$ on both sides, we obtain the
lemma.
\end{proof}

\begin{theorem} \label{basic}
Let $i>0$. If $\p_{0,0}^iz_{k,g}=0$ for $k\le K$, then $z_{k,g}=0$ for
$k\le K$.
\end{theorem}
\begin{proof}
We will show that if $\p_{0,0}^iz_{k,g}=0$ for all $k\le K$, then
$\p_{0,0}^{i-1}z_{k,g}=0$ for all $k\le K$. The theorem follows by downward
induction on $i$.

The puncture equation implies that $z_{-1,g}=0$ and hence that
$\p_{0,0}^{i-1}z_{-1,g}=0$, so by induction, we may suppose that
$\p_{0,0}^{i-1}z_{k-,g}=0$. Since $[\p_{0,0},\CL_{-1}]=0$,
Lemma~\ref{downward} implies that
$$
\CL_{-1}\p_{0,0}^{i-1}z_{k,g}=-(k+1)\p_{0,0}^{i-1}z_{k-1,g}=0 .
$$
Lemma \ref{key} now shows that $\p_{0,0}^{i-1}z_{k,g}$ is constant for all
$k\ge0$.

The dilaton operator $\Dil$ \eqref{Dilaton} commutes with $L_k$, and
$[\p_{0,0},\Dil]=\p_{0,0}$. This implies that
\begin{equation} \label{dil}
\Bigl( - \p_{1,0} + \sum_{m=0}^\infty t_m^a \p_{m,a} + (2g+i-3) \Bigr)
\p_{0,0}^{i-1} z_{k,g} = 0 .
\end{equation}
If $\p^{i-1}_{0,0}z_{k,g}$ is constant, we see that $(2g+i-3)
\p^{i-1}_{0,0} z_{k,g} = 0$, and hence, provided $2g+i\ne 3$, that
$\p^{i-1}_{0,0}z_{k,g}=0$ for all $k\ge0$.

There remain the exceptional cases $(g,i)=(0,3)$ and $(1,1)$. Suppose that
$\p_{0,0}^2z_{k,0}$ is constant for $k\ge0$. Then $\p_{0,0}^2z_{k,0}$ may
be calculated by applying the operator $\p_{0,0}^2$ to the explicit formula
\eqref{zkg} for $z_{k,0}$ and evaluating at $0\in\PHASE(V)$: if $k>r$, this
gives
\begin{multline*}
\p_{0,0}^2z_{k,0} = \sum_{i=0}^r \biggl( - [\tfrac{3-r}{2}]^k_i (R^i)^a_0
\<\tau_{k-i+1,a}\tau_{0,0}\tau_{0,0}\>_0^V \\ + 2 [\mu_a\!+\!\half]^k_i
(R^i)_0^a \<\tau_{k-i,a}\tau_{0,0}\>_0^V + \sum_{m=i-k}^{-1} (-1)^m
[\mu_a\!+\!m\!+\!\half]^k_i (R^i)^{ab} \\ \Bigl( \half
\<\tau_{-m-1,a}\tau_{0,0}\>_0^V \<\tau_{m+k-i,b}\tau_{0,0}\>_0^V +
\<\tau_{-m-1,a}\>_0^V \<\tau_{m+k-i,b}\tau_{0,0}\tau_{0,0}\>_0^V \Bigr)
\biggr) .
\end{multline*}
Choose $\beta\in H_2^+(V,\Z)$. By the dimension formula, the coefficient of
$q^\beta$ in each of the terms in the above formula vanishes unless
$r-1+c_1(V)\cap\beta=k$. It follows that for sufficiently large $k$, the
coefficient of $q^\beta$ in $\p_{0,0}^2z_{k,0}$ vanishes. By a downward
induction using Lemma \ref{downward}, it follows that the coefficient of
$q^\beta$ in $\p_{0,0}^2z_{k,0}$ vanishes for all $k$.

Similarly, suppose that $z_{k,1}$ is constant for $k\ge0$. Then evaluating
the explicit formula \eqref{zkg} for $z_{k,1}$ at $0\in\PHASE(V)$ gives
\begin{multline*}
z_{k,1} = \sum_{i=0}^{k+1} \biggl( - [\tfrac{3-r}{2}]^k_i (R^i)^b_0
\<\tau_{k-i+1,b}\>_1^V + \sum_{m=i-k}^{-1} (-1)^m
[\mu_a\!+\!m\!+\!\half]^k_i (R^i)^{ab} \\ \Bigl( \half
\<\tau_{-m-1,a}\tau_{m+k-i,b}\>_0^V + \<\tau_{-m-1,a}\>_0^V
\<\tau_{m+k-i,b}\>_1^V \Bigr) \biggr) .
\end{multline*}
By the dimension formula, the coefficient of $q^\beta$ in each of the terms
in the above formula vanishes unless $c_1(V)\cap\beta=k$. It follows that
for sufficiently large $k$, the coefficient of $q^\beta$ in $z_{k,1}$
vanishes, and we again conclude that $z_{k,1}=0$ for all $k$ by downward
induction using Lemma \ref{downward}.
\end{proof}

\section{The Virasoro conjecture in genus $0$} \label{genus0}

In this section, we present a proof of the Virasoro conjecture in genus
$0$. Our proof follows along the lines of the argument of Eguchi et al.\
\cite{EHX}; we have also borrowed some ingredients from the beautiful paper
of Dubrovin and Zhang \cite{DZ}, which proves a far-reaching generalization
of the genus $0$ Virasoro conjecture for any Frobenius manifold. (We
discuss some of their results in Section \ref{frobenius}.)

The first complete proof of the genus $0$ Virasoro conjecture of which we
are aware was given by Liu and Tian \cite{LiuT}; their proof of the
equation which Eguchi~et~al.\ call $\tilde{L}_1=0$ also influenced our
presentation.

Let $\CL_k$ be the vector field
\begin{equation} \label{CLk}
\CL_kf = \lim_{\hbar\to0} Z(V)^{-1} [ L_k , f ] Z(V) ,
\end{equation}
given by the explicit formula
\begin{align*}
\CL_k = \sum_{i=0}^{k+1} & \biggl( \sum_{m=i-k}^{-1} (-1)^m
[\mu_a\!+\!m\!+\!\half]^k_i (R^i)^{ab} \<\<\tau_{-m-1,a}\>\>_0^V
\p_{m+k-i,b} \\ & - [\tfrac{3-r}{2}]^k_i (R^i)^b_0 \p_{k-i+1,b} +
\sum_{m=0}^\infty [\mu_a\!+\!m\!+\!\half]^k_i (R^i)^b_a t^a_m \p_{m+k-i,b}
\biggr) .
\end{align*}
In particular, $\CL_{-1}$ is the vector field which we introduced in
Section \ref{Key}, while $\CL_0$ is given by the formula
$$
\CL_0 = - \half (3-r) \p_{1,0} +\sum_{m=0}^\infty
(\mu_a\!+\!m\!+\!\half) t_m^a \p_{m,a} - R_0^b \p_{0,b} +
\sum_{m=1}^\infty R_a^b t_m^a \p_{m-1,b} + \frac{1}{2\hbar} R_{ab}
t^a_0 t^b_0 .
$$
Starting from the explicit formula
$$
\p_{0,a} z_{k,0} = \CL_k \<\<\tau_{0,a}\>\>_0^V + \sum_{i=0}^k
[\mu_a\!+\!\half]^k_i (R^i)^b_a \<\<\tau_{k-i,b}\>\>_0^V + (R^{k+1})_{ab}
t^b_0 ,
$$
we will show that $\p_{0,0}z_{k,0}=0$; Theorem \ref{basic} then implies
that $z_{k,0}=0$ for $k\ge0$.

The arguments of this section work equally well when the matrix $\mu$
in the above constraints is replaced with the matrix $\mu_s$ of
\eqref{mus}. Thus, an analogue of the Virasoro conjecture holds in
genus $0$ for the Gromov-Witten invariants of any compact symplectic
manifold --- indeed, more generally, for any Frobenius manifold
\cite{DZ}.

Central to our argument is the Laurent series $\theta(\zeta) =
\theta^a(\zeta)\o x_a \in \CO_{\PHASE(V)}\o H(V)$, where
$$
\theta^a(\zeta) = \sum_{m=0}^\infty \zeta^{-m-1} \eta^{ab}
\<\<\tau_{m,b}\>\>_0^V + \sum_{m=0}^\infty (-\zeta)^m \tilde{t}_m^a .
$$

Let $\nabla:\CO_{\PHASE(V)}\to\CO_{\PHASE(V)}\o H(V)$ be the differential
operator
$$
\nabla f = \eta^{ab} \p_{0,a}f \o x_b .
$$
We also need the gradient $\Theta(\zeta)= \nabla\theta(\zeta) \in
\CO_{\PHASE(V)}\o\End(H(V))$, with coefficients
$$
\Theta_b^a(\zeta) = \delta_b^a + \sum_{m=0}^\infty \zeta^{-m-1} \eta^{ac}
\<\<\tau_{m,c}\tau_{0,b}\>\>_0^V .
$$

Observe that the equation $\zeta\CL_{-1}\theta(\zeta)+\theta(\zeta)=0$
holds. Since $[\nabla,\CL_{-1}]=0$, we also see that
$\zeta\CL_{-1}\Theta(\zeta)+\Theta(\zeta)=0$.

Denote the matrix $\Res_\zeta(\Theta)\in\CO_{\PHASE(V)}\o\End(H(V))$
by $\CC$; it has coefficients
$$
\CC_b^a = \eta^{ac}\<\<\tau_{0,c}\tau_{0,b}\>\>_0^V .
$$
Denote the matrix $-\CL_0\CC$ by $\DC$.

A basic property of Gromov-Witten invariants in genus $0$ is the
topological recursion relation
\begin{equation} \label{trr0}
\<\<\tau_{k,a}\tau_{\ell,b}\tau_{m,c}\>\>_0^V = \eta^{ef}
\<\<\tau_{k,a}\tau_{0,e}\>\>_0^V
\<\<\tau_{0,f}\tau_{\ell,b}\tau_{m,c}\>\>_0^V ;
\end{equation}
this is ultimately a consequence of the fact that the tautological line
bundles $\Om_i$ vanish on the zero dimensional variety $\Mbar_{0,3}$. This
relation has two consequences ((6.28) and (6.31) of Dubrovin \cite{D}). The
first of these is called the \emph{quantum differential equation} by
Givental.

\begin{lemma} \label{quantum}
If $\xi$ is a vector field on the large phase space,
$$
\zeta \CL_\xi \Theta(\zeta) = \Theta(\zeta) \CL_\xi\CC .
$$
\end{lemma}

\begin{lemma} \label{invert}
$\Theta(\zeta)\Theta^*(-\zeta)=\Id$
\end{lemma}
\begin{proof}
This follows by induction from the formula
$$
\eta^{ef} \<\<\tau_{k,a}\tau_{0,e}\>\>_0^V
\<\<\tau_{0,f}\tau_{\ell,b}\>\>_0^V =
\<\<\tau_{k,a}\tau_{\ell+1,b}\>\>_0^V +
\<\<\tau_{k+1,a}\tau_{\ell,b}\>\>_0^V .
$$
This formula holds at $0\in\PHASE(V)$, since both sides vanish there, while
by \eqref{trr0},
$$
\p_{m,c} \Bigl( \eta^{ef} \<\<\tau_{k,a}\tau_{0,e}\>\>_0^V
\<\<\tau_{0,f}\tau_{\ell,b}\>\>_0^V - \<\<\tau_{k,a}\tau_{\ell+1,b}\>\>_0^V
- \<\<\tau_{k+1,a}\tau_{\ell,b}\>\>_0^V \Bigr) = 0
$$
for all $m\ge0$ and $c\in A$.
\end{proof}

Let $\UU(\zeta)$ be the generating function $\UU(\zeta) = \theta^*(-\zeta)
\Theta(\zeta)$, and let $\UU[n]$ be the coefficient of $\zeta^n$ in $\UU$.
These elements of $\CO_{\PHASE(V)}\o H(V)$ were introduced by Eguchi,
Yamada and Yang \cite{EYY} in their study of Gromov-Witten invariants in
higher genus.

The main result of this section is that $\UU[n]=0$ for $n\le0$. For $n=0$,
this follows from the puncture equation, since $\UU[0]=\nabla
z_{-1,0}$. The equation $\UU[-1]=0$ is a consequence of the dilaton
equation, since
$$
\UU[-1] = \nabla \Bigl( \sum_{m=0}^\infty \tilde{t}_m^a
\<\<\tau_{m,a}\>\>_0^V - 2 \<\<~\>\>_0^V \Bigr) = 0 .
$$
The equation $\UU[-2]=0$ is a consequence of the equation $\tilde{L}_1=0$
of Eguchi et al.\ \cite{EHX}, for which a proof has been given by Liu and
Tian \cite{LiuT}:
$$
\UU[-2] = \nabla \Bigl( \sum_{m=0}^\infty \tilde{t}_m^a
\<\<\tau_{m+1,a}\>\>_0^V - \frac12 \eta^{ab} \<\<\tau_{0,a}\>\>_0^V
\<\<\tau_{0,b}\>\>_0^V \Bigr) = 0 .
$$
More generally, using the Grothendieck-Riemann-Roch theorem, Faber and
Pandharipande \cite{FP} prove the equations
\begin{equation} \label{FP}
\sum_{m=0}^\infty \tilde{t}_m^a \<\<\tau_{m+2\ell-1,a}\>\>_0^V + \frac12
\sum_{m=1-2\ell}^{-1} (-1)^m \eta^{ab} \<\<\tau_{-m-1,a}\>\>_0^V
\<\<\tau_{m+2\ell-1,b}\>\>_0^V = 0 .
\end{equation}
Applying $\nabla$ yields the equations $\UU[-2\ell]=0$.

However, for $n<-1$ odd, the function $\UU[-n]$ does not have the form
$\nabla z$; in these cases, the following result is new.
\begin{proposition} \label{P}
For $n\le0$, $\UU[n]$ vanishes.
\end{proposition}
\begin{proof}
As we mentioned already, the constant term $\UU[0]$ of $\UU(\zeta)$ equals
$\nabla z_{-1,0}$, and vanishes by the puncture equation. We now argue by
downward induction on $n$. We have
$$
\nabla \UU(\zeta) = \nabla\theta^*(-\zeta) \Theta(\zeta) + \theta^*(-\zeta)
\nabla \Theta(\zeta) = \Theta^*(-\zeta)\Theta(\zeta) + \zeta^{-1}
\UU(\zeta) \nabla\CC .
$$
Since by Lemma \ref{invert}, $\Theta^*(-\zeta)\Theta(\zeta)=\Id$, we see
that $\nabla \UU[-n] = \UU[1-n]\nabla\CC$.

Since $\zeta\CL_{-1}\theta(\zeta)+\theta(\zeta)=0$ and
$\zeta\CL_{-1}\Theta(\zeta)+\Theta(\zeta)=0$, we see that
$\CL_{-1}\UU(\zeta)=0$. Lemma \ref{key} shows that the vanishing of
$\UU[1-n]$ implies that $\UU[-n]$ is constant. But $\theta(\zeta)$ vanishes
at $0\in\PHASE(V)$, showing that $\UU[-n]$ does too.
\end{proof}

Introduce the differential operator on the circle:
$$
\delta = - \zeta^2 \p + \zeta (\mu-\half) + R .
$$
It is straightfoward, if a little tedious, to show that $\nabla z_{k,0}$ is
the constant term of $\theta^*(-\zeta) \delta^{k+1} \Theta(\zeta)$.  The
operator $\delta$ was introduced in Eguchi, Hori and Xiong \cite{EHX0},
where the following formula is proved. (Recall that $\DC=-\CL_0\CC$.)
\begin{proposition} \label{UV}
$\Theta(\zeta)^{-1} \delta\Theta(\zeta) = \zeta (\mu-\half) + \DC$
\end{proposition}
\begin{proof}
On the one hand, Lemma \ref{quantum} implies that
$\zeta\CL_0\Theta(\zeta)=\Theta(\zeta)\CL_0\CC$. On the other hand, Hori's
equation $z_{0,0}=0$ implies that
\begin{align*}
0 &= \p_{n,a} \p_{0,b} z_{0,0} \\ &= \begin{cases} 
\bigl( \CL_0 + \mu_a +\mu_b + 1 \bigr) \<\<\tau_{n,a}\tau_{0,b}\>\>_0^V +
R_{ab} , & n=0 , \\[3pt]
\bigl( \CL_0 + n + \mu_a + \mu_b + 1 \bigr)
\<\<\tau_{n,a}\tau_{0,b}\>\>_0^V + R_a^e \<\<\tau_{n-1,e}\tau_{0,b}\>\>_0^V
, & n>0 ;
\end{cases}
\end{align*}
in other words, $\CL_0 \Theta(\zeta) + \zeta^{-1} \delta\Theta(\zeta) -
\Theta(\zeta) (\mu-\half)=0$.
\end{proof}

The explicit formula
\begin{equation} \label{DC}
\DC = \CC + R + [\mu,\CC]
\end{equation}
follows by taking the constant term of Proposition \ref{UV}.

It is now easy to see that $\nabla z_{k,0}$, and hence $\p_{0,0}z_{k,0}$,
vanishes for $k\ge0$. Iterating Proposition \ref{UV}, we see that there are
matrices $\QQ_{k,i}$ of functions on $\PHASE(V)$ such that
\begin{equation} \label{iterate}
\delta^{k+1}\Theta(\zeta) = \sum_{i=0}^{k+1} \zeta^i \Theta(\zeta)
\QQ_{k,i} .
\end{equation}
We conclude that $\nabla z_{k,0} = \sum_{i=0}^{k+1} \UU[-i] \QQ_{k,i}$,
which vanishes by Proposition \ref{P}.

The matrices $\QQ_{k,i}$ may be calculated recursively, starting with
$\QQ_{-1,i}=\delta_{i,0}$:
$$
\QQ_{k,i} = (\mu+\half-i) \QQ_{k-1,i-1} + \DC \QQ_{k-1,i} .
$$
In particular, $\QQ_{k,0}=\DC^{k+1}$.

\section{Frobenius manifolds and the Virasoro conjecture} \label{frobenius}

Dubrovin introduced Frobenius manifolds as an axiomatization of the
structure of Gromov-Witten invariants in genus $0$. Dubrovin and Zhang have
shown that the Virasoro conjecture in genus $0$ has a reformulation in the
language of Frobenius manifolds; the methods developed in the last section
allow an efficient proof of this relationship.

Let $\phase$ be a smooth (super)scheme (or manifold) with structure sheaf
$\CO_\phase$ and tangent sheaf $\CT_\phase$. A pre-Frobenius structure on
$\phase$ consists of the following data:
\begin{enumerate}
\item a (graded) commutative product $\CT_\phase\o \CT_\phase\to
\CT_\phase$, which we denote by $X\o Y\mapsto X\circ Y$;
\item a non-degenerate (graded) symmetric bilinear form
$\CT_\phase\o\CT_\phase\to\CO_\phase$ (i.e. pre-Riemannian metric), which
we denote by $X\o Y\mapsto (X,Y)$, compatible with the product in the sense
that $(X,Y\circ Z) = (X\circ Y,Z)$;
\item An Euler vector field $\Euler$, that is, a linear vector field
$\nabla\nabla\Euler=0$ which defines a grading for the product,
$$
[\Euler , X\circ Y ] = [\Euler, X ]  \circ Y + X \circ [ \Euler , Y ] +
X\circ Y ,
$$
and which is conformal: there is a constant $r$ such that
$$
\Euler (X,Y) = ( [\Euler,X] , Y ) + ( X , [\Euler,Y] ) + (2-r)(X,Y) .
$$
\end{enumerate}

Introduce the pencil of connections $\nabla^\lambda_XY = \nabla_XY +
\lambda X\circ Y$, where $\nabla$ is the Levi-Civita connection associated
to the bilinear form $(X,Y)$, that is, the unique torsion-free connection
such that
$$
\nabla_Z (X,Y) = ( \nabla_Z X , Y ) + ( X , \nabla_Z,Y ) .
$$
\begin{definition} \label{Frobenius}
A \emph{Frobenius manifold} $\phase$ is a manifold with pre-Frobenius
structure such that $\nabla^\lambda$ is flat for all $\lambda$.
\end{definition}
Definition \ref{Frobenius} amounts to the conditions that the Levi-Civita
connection $\nabla$ is flat and that there is a function $\Phi$ on $\phase$
such that the symmetric three-tensor $(X\circ Y,Z)$ is the third derivative
of $\Phi$.

The basic example of a Frobenius manifold is the small phase space
$\phase(V)$ of a smooth projective variety $V$. Recall the
Witten-Dijkgraaf-Verlinde-Verlinde (WDVV) equation.
\begin{proposition}
$$
\eta^{ef} \<\<\tau_{k,a}\tau_{\ell,b}\tau_{0,e}\>\>_0^V
\<\<\tau_{0,f}\tau_{m,c}\tau_{n,d}\>\>_0^V = \eta^{ef}
\<\<\tau_{k,a}\tau_{m,c}\tau_{0,e}\>\>_0^V
\<\<\tau_{0,f}\tau_{\ell,b}\tau_{n,d}\>\>_0^V
$$
\end{proposition}
\begin{proof}
Apply the vector field $\p_{n,d}$ to both sides of the genus $0$
topological recursion relation \eqref{trr0}:
\begin{align*}
\<\<\tau_{k+1,a}\tau_{\ell,b}\tau_{m,c}\tau_{n,d}\>\>_0^V &= \eta^{ef}
\<\<\tau_{k,a}\tau_{\ell,b}\tau_{0,e}\>\>_0^V
\<\<\tau_{0,f}\tau_{m,c}\tau_{n,d}\>\>_0^V \\ &+ \eta^{ef}
\<\<\tau_{k,a}\tau_{\ell,b}\tau_{n,d}\tau_{0,e}\>\>_0^V
\<\<\tau_{0,f}\tau_{m,c}\>\>_0^V .
\end{align*}
The left-hand side and the first term of the right-hand side are invariant
under exchange of $\tau_{\ell,b}$ and $\tau_{m,c}$, from which the result
follows.
\end{proof}

(We have stated the WDVV equation in the case where there are no odd
cohomology classes; in general, we must multiply by a sign determined in
the usual way by the sign convention for $\Z/2$-graded vector spaces.)

We may now define a Frobenius structure on $\phase(V)$. The metric on
$\phase(V)$ is the flat metric $(\p_a,\p_b)=\eta_{ab}$. The small phase
space $\phase(V)$ may be embedded into the large phase space $\PHASE(V)$
along the subscheme $\{t_m^a=0\mid m>0\}$, by sending the point with
coordinates $u^a$ to the point with coordinates $t_m^a=\delta_{m,0}u^a$;
call this embedding $s$. The product on $\CT_{\phase(V)}$ is given by the
formula
$$
\p_a \circ \p_b = \eta^{ef} s^*\<\<\tau_{0,a}\tau_{0,b}\tau_{0,e}\>\>_0^V
\p_f .
$$
The function $\Phi=s^*\<\<~\>\>_0^V$ is a potential for this product. To
complete the construction of the Frobenius manifold, it only remains to
construct the Euler vector field.
\begin{proposition} \label{Euler}
The vector field $\CE = (1-p_a) u^a \p_a + R_0^a \p_a$ is an Euler vector
field on $\phase(V)$, with $r=\dim_\C(V)$.
\end{proposition}
\begin{proof}
Since $[\CE,\p_a]=(p_a-1)\p_a$, we see that
$$
\Euler (\p_a,\p_b) - ( [\Euler,\p_a] , \p_b ) - ( \p_a , [\Euler,\p_b] ) =
\bigl( (1-p_a) + (1-p_b) \bigr) \eta_{ab} = (2-r) \eta_{ab} .
$$

Let $\CL_0$ be the vector field $L_0-\rho(V)$ on $\PHASE(V)$. The equations
$\CL_0\<\<~\>\>_0^V=0$ and $\Dil\<\<~\>\>_0^V=2\<\<~\>\>_0^V$ imply that
\begin{equation} \label{eul}
(\CL_0+\half(r-3)\Dil)\<\<\tau_{0,a}\tau_{0,b}\tau_{0,c}\>\>_0^V +
(p_a+p_b+p_c-r) \<\<\tau_{0,a}\tau_{0,b}\tau_{0,c}\>\>_0^V = 0 .
\end{equation}
The vector field $\CL_0+\half(r-3)\Dil$ is tangential to the image of the
emebedding $s$; pulling the identity \eqref{eul} back to $\phase(V)$ by
$s$, we see that
\begin{equation} \label{euler}
\bigl( - \Euler + p_a + p_b + p_c - r \bigr) ( \p_a\circ\p_b , \p_c ) = 0 .
\end{equation}
On the other hand,
\begin{multline*}
\bigl( [\Euler , \p_a\circ \p_b ] - [\Euler , \p_a ] \circ \p_b - \p_a
\circ [ \Euler , \p_b ] , \p_c \bigr) \\ = \bigr( \Euler + (r-p_c-1) -
(p_a-1) - (p_b-1) \bigr) ( \p_a\circ\p_b , \p_c ) .
\end{multline*}
By \eqref{euler}, this equals $( \p_a\circ\p_b , \p_c )$.
\end{proof}

There is a fibration $u$ of the large phase space $\PHASE(V)$ over the
small phase space $\phase(V)$, obtained by mapping the point with
coordinates $t_m^a$ to the point with coordinates
$u^a=\eta^{ab}\<\<\tau_{0,b}\>\>_0^V=t_0^a+O(|t|^2)$. The importance
of the fibration $u:\PHASE(V)\to\phase(V)$ is illustrated by the
formula $\Theta(\zeta)=u^*s^*\Theta(\zeta)$ of Dijkgraaf and Witten
\cite{DW} (cf.\ Section 7 of \cite{genus2}).
\begin{definition}
A vector field $\xi$ on $\PHASE(V)$ is \emph{horizontal} if
$$
\CL_\xi u^*\CO_{\phase(V)}\subset u^*\CO_{\phase(V)} ;
$$
the vector field thereby induced on $\phase(V)$ is denoted $u_*\xi$.
\end{definition}

\begin{lemma} \label{horizontal}
If $\xi$ is a horizontal vector field on the large phase space,
$$
u_*\xi\circ X = \bigl(s^*\CL_\xi\CC\bigr)X .
$$
\end{lemma}
\begin{proof}
Since $u$ is a submersion, we may add a vertical vector field to $\xi$ so
that it is tangential to the section $s$ of $u$. (For example, this is what
we did when we replaced $\CL_0$ by $\CL_0+\half(r-3)\Dil$ in the proof of
Proposition \ref{Euler}.) If $\xi$ satisfies this additional condition,
then
$$
\xi = \sum_{m=0}^\infty \xi_m^a \p_{m,a}
$$
where $s^*\xi_m^a=0$ if $m>0$. It follows that
$$
s^*\CL_\xi\CC_{ab} = s^* \bigl( \xi_0^f
\<\<\tau_{0,f}\tau_{0,a}\tau_{0,b}\>\>_0^V \bigr) = ( u_*\xi\circ\p_a ,
\p_b ) . \qed
$$
\def\qed{}
\end{proof}

The Virasoro conjecture in genus $0$ implies that the vector fields $\CL_k$
on the large phase space $\PHASE(V)$ introduced in \eqref{CLk} are given by
the formulas
$$
\CL_k = \lim_{\hbar\to0} Z(V)^{-1} \circ L_k \circ Z(V) .
$$
This implies the following lemma.
\begin{lemma}
The vector fields $\CL_k$ satisfy the Virasoro commutation relations
$$
[\CL_k,\CL_\ell] = (k-\ell) \CL_{k+\ell} .
$$
\end{lemma}

We can now prove the important result of Dubrovin and Zhang \cite{DZ},
which relates the Virasoro conjecture in genus $0$ to the Frobenius
geometry of the small phase space.
\begin{theorem} \label{DZ}
The vector fields $\CL_k$ on the large phase space are horizontal, and
$u_*\CL_k + \CE^{\circ(k+1)}=0$. Equivalently,
$\zeta\CL_k\Theta(\zeta)+\Theta(\zeta)\DC^{k+1}=0$.
\end{theorem}
\begin{proof}
We must show that $\CL_ku^a$ depends only on the functions $u^b$. We start
with $k=-1$. The coefficient of $\zeta$ in the equation
$\zeta\CL_{-1}\Theta(\zeta)+\Theta(\zeta)=0$ is the equation
$$
\CL_{-1}\CC + \Id = 0 .
$$
The first row of this formula is $\CL_{-1}u^a + \delta_0^a = 0$, which
shows that $\CL_{-1}$ is horizontal, and that $u_*\CL_{-1}+\p_0=0$. But the
vector field $\p_0$ is the identity vector field on the Frobenius manifold
$\phase(V)$: $\p_0\circ X=X$.

We next turn to $k=0$. The first row of the formula \eqref{DC} is
$$
\CL_0 u^a + (1-p_a) u^a + R_0^a = 0 .
$$
This shows that $\CL_0$ is horizontal, that $u_*\CL_0 + \Euler = 0$, and,
applying Lemma \ref{horizontal}, that $\Euler \circ X=s^*\DC X$.

Granted the equation $z_{k,0}=0$, it is straightfoward to show that
$-\CL_k\CC$ equals the coefficient of $\zeta^0$ in the generating function
$\Theta^*(-\zeta) \delta^{k+1} \Theta(\zeta)$. By \eqref{iterate}, it
follows that $\CL_k\CC=-\QQ_{k,0}=-\DC^{k+1}$, and hence that
$$
\Euler^{\circ(k+1)} \circ X = s^*\DC^{k+1} X = - \bigl( s^*\CL_k\CC \bigr)
X .
$$
Taking $X=\p_0$, the theorem follows.
\end{proof}

It follows from this theorem that the vector fields $\Euler^{\circ k}$ on
$\phase(V)$ satisfy the Virasoro commutation relations
$$
\bigl[ \Euler^{\circ k} , \Euler^{\circ\ell} \bigr] = (\ell-k)
\Euler^{\circ(k+\ell-1)} .
$$
These commutation relations have recently been proved for all Frobenius
manifolds by Hertling and Manin \cite{HM}.

\section{The Virasoro conjecture for a point and the KdV hierarchy}
\label{point}

If $V$ is a point, $\Mbar_{g,n}(V,0)$ is the moduli space $\Mbar_{g,n}$ of
$n$-pointed stable curves of arithmetic genus $g$ introduced by Deligne,
Mumford and Knudsen. In this case, the moduli stack is smooth for all $g$
and $n$, and has dimension equal to its virtual dimension, namely
$\dim\Mbar_{g,n}=3(g-1)+n$. In particular, the virtual fundamental class is
just the fundamental class $[\Mbar_{g,n}]\in
H_{6(g-1)+2n}(\Mbar_{g,n},\Q)$, and the Gromov-Witten invariants are just
the intersection numbers
$$
\< \tau_{k_1} \dots \tau_{k_n} \>_g = \int_{\Mbar_{g,n}} \Psi_1^{k_1} \dots
\Psi_n^{k_n} .
$$

For $\Mbar_{0,3}$ and $\Mbar_{1,1}$, it is easy to calculate these
intersection numbers directly.

\smallskip

{$\bullet$}
The stack $\Mbar_{0,3}$ has dimension zero, so the only
intersection number to be calculated is $\<\tau_0^3\>_0$; since
$\Mbar_{0,3}$ consists of a single point, we see that $\<\tau_0^3\>_0=1$.

\smallskip

{$\bullet$}
The stack $\Mbar_{1,1}$ has dimension one, so the only
intersection number to be calculated is $\<\tau_1\>_1$. There are many ways
to do this: for example, we may identify $\Mbar_{1,1}$ with the
compactification of the moduli space of elliptic curves by a single cusp
and sections of $\Om_1^n$ over $\Mbar_{1,1}$ with cusp forms of weight $n$
for the modular group $\SL(2,\Z)$. The cusp form
$$
\Delta = q \prod_{n=1}^\infty (1-q^n)^{24} , \quad \text{where $q=e^{2\pi
i\tau}$,}
$$
of weight $12$ is nonzero everywhere except at the cusp, where it has a
simple zero. The stable curve represented by the cusp has a non-trivial
involution, and hence the associated divisor has degree $\half$; this shows
that the line bundle $\Om_1^{12}$ on $\Mbar_{1,1}$ has degree $\half$, and
hence that $\<\tau_1\>_1 = \tfrac{1}{24}$.

\medskip

When $V$ is a point, the large phase space $\PHASE$ has coordinates
$\{t_m\mid m\ge0\}$; denote the Gromov-Witten potential in this case by
$Z$. The puncture and dilaton equations amount to the following:
\begin{align*}
\<\<\tau_0\>\>_g &= \sum_{m=1}^\infty t_m \<\<\tau_{m-1}\>\>_g +
\delta_{g,0} \frac{t_0^2}{2} , \\
\<\<\tau_1\>\>_g &= \sum_{m=0}^\infty t_m \<\<\tau_m\>\>_g +
\frac{\delta_{g,1}}{24} .
\end{align*}
It follows from the puncture equation that
$$
\<\tau_{k_1}\dots\tau_{k_n}\>_0 = \begin{cases} \dfrac{(n-3)!}{k_1!\dots
k_n!} , & k_1+\dots+k_n=n-3 , \\[10pt] 0 , & \text{otherwise.}
\end{cases}
$$
The puncture and dilaton equations together allow us to express all of the
Gromov-Witten invariants $\<\tau_{k_1}\dots\tau_{k_n}\>_g$ in terms of
those with $k_i>1$: in genus $1$,
$$
\<\<~\>\>_1 = \frac{1}{24} \log u' ,
$$
where $u'=\<\<\tau_0{}^3\>\>_0$, and in higher genus,
\begin{equation} \label{IZ}
\<\<~\>\>_g = \sum_{n=1}^{3g-3} \frac{1}{n!} \bigl(u'\bigr)^{-(2g-2+n)}
\sum_{\substack{k_1+\dots+k_n=3g-3+n\\k_i>1}}
\<\tau_{k_1}\dots\tau_{k_n}\>_g \UU[k_1] \dots \UU[k_n] ,
\end{equation}
where, as in Section 4,
$$
\UU[k] = t_k + \sum_{m=0}^\infty t_{m+k+1} \<\<\tau_0\tau_m\>\>_0 .
$$
(For the proof of this formula, see Section 5 of Itzykson and Zuber
\cite{IZ}.) In particular, the Gromov-Witten invariants in genus $g$ are
determined by $p(3g-3)$ numbers, where $p(3g-3)$ is the number of
partitions of $3g-3$.

\subsection{Calculation of $\<\<~\>\>_2$}

The equations $z_{k,g}=0$, or equivalently,
\begin{multline*}
\tfrac{\Gamma(k+\frac52)}{\Gamma(\frac32)} \<\<\tau_{k+1}\>\>_g =
\sum_{m=0}^\infty \tfrac{\Gamma(m+k+\frac32)}{\Gamma(m+\frac12)} t_m
\<\<\tau_{m+k}\>\>_g \\ + \tfrac{1}{2} \sum_{m=-k}^{-1} (-1)^m
\tfrac{\Gamma(m+k+\frac32)}{\Gamma(m+\frac12)} \Bigl(
\<\<\tau_{-m-1}\tau_{m+k}\>\>_{g-1} + \sum_{h=0}^g \<\<\tau_{-m-1}\>\>_h
\<\<\tau_{m+k}\>\>_{g-h} \Bigr) ,
\end{multline*}
may be used to inductively determine all of the intersection numbers
$\<\tau_{k_1}\dots\tau_{k_n}\>_g$. Let us illustrate how this scheme works
in genus $2$, by calculating the intersection numbers $\<\tau_4\>_2$,
$\<\tau_2\tau_3\>_2$ and $\<\tau_2{}^3\>_2$.

The equation $z_{3,2}=0$ gives
\begin{align*}
\tfrac{945}{16} \<\<\tau_4\>\>_2 &= \sum_{m=0}^\infty
\tfrac{\Gamma(m+\frac92)}{\Gamma(m+\frac12)} t_m \<\<\tau_{m+3}\>\>_2 \\
& \quad + \tfrac{15}{16} \bigl( \<\<\tau_0\>\>_0\<\<\tau_2\>\>_2
+ \<\<\tau_0\>\>_1\<\<\tau_2\>\>_1 + \<\<\tau_0\>\>_2\<\<\tau_2\>\>_0 +
\<\<\tau_0\tau_2\>\>_1 \bigr) \\
& \quad + \tfrac{9}{32} \bigl( 2\<\<\tau_1\>\>_0\<\<\tau_1\>\>_2
+ \<\<\tau_1\>\>_1\<\<\tau_1\>\>_1 + \<\<\tau_1{}^2\>\>_1 \bigr)
\end{align*}
Setting the variables $t_m$ to zero, we see that
$$
\tfrac{945}{16} \<\tau_4\>_2 = \tfrac{15}{16} \<\tau_0\tau_2\>_1 +
\tfrac{9}{32} \bigl( \<\tau_1\>_1 \<\tau_1\>_1 + \<\tau_1{}^2\>_1 \bigr) =
\tfrac{105}{2048} ,
$$
and hence that $\<\tau_4\>_2=\frac{1}{1152}$.

The equation $z_{2,2}=0$ gives
\begin{align*}
\tfrac{105}{8} \<\<\tau_3\>\>_2 &= \sum_{m=0}^\infty
\tfrac{\Gamma(m+\frac72)}{\Gamma(m+\frac12)} t_m \<\<\tau_{m+2}\>\>_2 \\ &
\quad + \tfrac{3}{8} \bigl( \<\<\tau_0\>\>_0 \<\<\tau_1\>\>_2 +
\<\<\tau_0\>\>_1 \<\<\tau_1\>\>_1 + \<\<\tau_1\>\>_0 \<\<\tau_0\>\>_2 +
\<\<\tau_0\tau_1\>\>_1 \bigr) .
\end{align*}
Applying the operator $\p_2$ and setting all of the variables $t_m$ to
zero, we see that
$$
\tfrac{105}{8} \<\tau_2\tau_3\>_2 = \tfrac{315}{8} \<\tau_4\>_2 +
\tfrac{3}{8} \bigl( \<\tau_0\tau_2\>_1 \<\tau_1\>_1 +
\<\tau_0\tau_1\tau_2\>_1 \bigr) = \tfrac{203}{3072} ,
$$
so that $\<\tau_2\tau_3\>_2=\frac{29}{5760}$.

Finally, the equation $z_{1,2}=0$ gives
$$
\tfrac{15}{4} \<\<\tau_2\>\>_2 = \sum_{m=0}^\infty
\tfrac{\Gamma(m+\frac52)}{\Gamma(m+\frac12)} t_m \<\<\tau_{m+1}\>\>_2 +
\tfrac{1}{8} \bigl( 2 \<\<\tau_0\>\>_0 \<\<\tau_0\>\>_2 + \<\<\tau_0\>\>_1
\<\<\tau_0\>\>_1 + \<\<\tau_0{}^2\>\>_1 \bigr) .
$$
Applying the operator $\p_2{}^2$ and setting all of the variables
$t_m$ to zero, we see that
$$
\tfrac{15}{4} \<\tau_2{}^3\>_2 = \tfrac{35}{2} \<\tau_2\tau_3\>_2 +
\tfrac{1}{4} \<\tau_0\tau_2\>_1 \<\tau_0\tau_2\>_1 + \tfrac{1}{8}
\<\tau_0{}^2\tau_2{}^2\>_1 = \tfrac{7}{64} ,
$$
so that $\<\tau_2{}^3\>_2=\frac{7}{240}$. These results agree with the
calculations of Mumford \cite{Mumford}.

By \eqref{IZ}, we conclude that
$$
\<\<~\>\>_2 = \frac{1}{1152} \frac{\UU[4]}{(u')^3} +
\frac{29}{5760} \frac{\UU[2]\UU[3]}{(u')^4} +
\frac{7}{1440} \frac{\UU[2]^3}{(u')^5} .
$$

\subsection{Gelfand-Dikii polynomials}
We now turn to the equivalence discovered by Dijkgraaf, Verlinde and
Verlinde \cite{DVV} between the Witten and Virasoro conjectures for the
Gromov-Witten invariants of a point. To state the Witten conjecture, we
first recall the definition of the Gelfand-Dikii polynomials, introduced in
\cite{GD}. These are a sequence of differential polynomials
$$
R_m(\u) \in \Q_\hbar\{\u\} = \Q[\hbar][\u^{(i)}\mid i\ge0]
$$
associated to the asymptotic expansion for small time of the heat-kernel of
a Sturm-Liouville operator.

Let $\p$ be the derivation on $\Q_\hbar\{\u\}$, defined on the generators
by $\p\u^{(i)}=\u^{(i+1)}$.
\begin{lemma}
If $f$ satisfies the Sturm-Liouville equation
$$
\bigl( \tfrac{\hbar}{2} \p^2 + \u \bigr) f = zf ,
$$
then $\HH(f^2) = z\p(f^2)$, where $\HH$ is the third-order linear
differential operator
$$
\HH = \tfrac{\hbar}{8} \p^3 + \u\p + \half \u' .
$$
\end{lemma}

The differential polynomials $R_m(\u)$, $m>0$, are defined by the recursion
\begin{equation} \label{GD}
\HH R_m = \bigl( m+\half \bigr) \p R_{m+1} ,
\end{equation}
where $R_0(\u)=1$, while the constant term of $R_m(\u)$ vanishes for $m>0$.
For example,
\begin{align*}
R_1 &= \u , \\
R_2 &= \tfrac{\hbar}{12}\u^{(2)} + \tfrac{1}{2}\u^2 , \\
R_3 &= \tfrac{\hbar^2}{240}\u^{(4)} + \tfrac{\hbar}{12}\u\u^{(2)} +
\tfrac{\hbar}{24}\bigl(\u'\bigr)^2 + \tfrac{1}{6}\u^3 .
\end{align*}
Note that $R_m(\u)$ is independent of $\u^{(i)}$ if $i>2m-2$.

\subsection{Witten's conjecture}
The function $\u=\hbar\<\<\tau_0{}^2\>\>$ on the large phase-space is a
``quantization'' of the function $u=\<\<\tau_0\tau_0\>\>_0$ which arose in
the study of Gromov-Witten invariants in genus $0$. Identify the
differential $\p$ of the algebra of differential polynomials
$\Q_\hbar\{\u\}$ with the differential $\p_0$ on the large phase space, so
that $\u^{(i)}=\p_0^i\u=\hbar\<\<\tau_0{}^{i+2}\>\>$. Witten's conjecture
asserts that the Gromov-Witten invariants of a point satisfy the equations
\begin{equation} \label{witten-conjecture}
\hbar \<\<\tau_m\tau_0\>\> = R_{m+1}(\u) .
\end{equation}
Applying the derivative $\p$, these equations may be written $\p_m \u = \p
R_{m+1}(\u)$. These are the equations of the KdV hierarchy; the first few
are
\begin{align*}
\p_0 \u &= \u' , \\
\p_1 \u &= \tfrac{\hbar}{12} \u^{(3)} + \u\u' , \\
\p_2 \u &= \tfrac{\hbar^2}{240} \u^{(5)} + \tfrac{\hbar}{6} \u' \u^{(2)} +
\tfrac{\hbar}{12} \u\u^{(3)} + \tfrac{1}{2} \u^2\u' .
\end{align*}
In combination with the puncture equation, this conjecture suffices to
determine the Gromov-Witten potential $Z$.

Witten's conjecture was proved by Kontsevich \cite{K}; see Itzykson and
Zuber \cite{IZ} and Looijenga \cite{L} for illuminating discussions of the
proof. Let $\CV_N$ be the space of $N\times N$ Hermitian matrices, and
given a positive-definite Hermitian matrix $\Lambda$, let $d\mu_\Lambda$ be
the probability measure on $\CV_N$ with density
$$
d\mu_\Lambda = \frac{1}{c_\Lambda} \exp\bigl(-\tfrac12\Tr(\Lambda M^2)\bigr)
.
$$
Kontsevich shows that the matrix integral
$$
Z_N(\Lambda) = \int_{\CV_N} \exp\bigl(\tfrac{i}{6}\Tr(M^3) \bigr)
d\mu_\Lambda
$$
depends on $\Lambda$ only through the variables
$$
t_m = - (2m-1)!! \Tr\bigl( \Lambda^{-2m-1} \bigr) , \quad m<N/2 ,
$$
and that
$$
\lim_{N\to\infty} Z_N(t_m) = Z .
$$
Using this representation, he shows that $Z$ satisfies the KdV hierarchy
\eqref{witten-conjecture}, thus proving Witten's conjecture.

In studying the Gromov-Witten invariants of a point, it is convenient
to employ rescaled coordinates on the large phase space:
$$
s_m = \frac{\Gamma(\frac{3}{2})}{\Gamma(m+\frac{3}{2})} t_m \quad,\quad
\tilde{s}_m = s_m - \tfrac{2}{3} \delta_{m,1} =
\frac{\Gamma(\frac{3}{2})}{\Gamma(m+\frac{3}{2})} \tilde{t}_m .
$$
The corresponding partial derivatives of the total potential are
$$
\<\<\sigma_{k_1}\dots\sigma_{k_n}\>\> = \frac{\p^n\log Z}{\p s_{k_1}\dots\p
s_{k_n}} .
$$
In this coordinate system, Witten's conjecture \eqref{witten-conjecture}
becomes the recursion
\begin{equation} \label{f}
\HH \p \<\<\sigma_{k-1}\>\> = \p^2 \<\<\sigma_k\>\> , \quad k>0 .
\end{equation}

\subsection{The Virasoro constraints}
Let $z_k$ be the Virasoro constraint $Z^{-1}L_kZ$. In terms of the
variables $s_m$, the Virasoro constraints with $k\ge0$ have the explicit
formulas
\begin{align*}
z_k &= - \<\<\sigma_{k+1}\>\> + \sum_{m=0}^\infty (m+\half) s_m
\<\<\sigma_{m+k}\>\> + \tfrac\hbar8 \sum_{i+j=k-1} \Bigl(
\<\<\sigma_i\sigma_j\>\> + \<\<\sigma_i\>\> \<\<\sigma_j\>\> \Bigr) \\
&= \sum_{m=0}^\infty (m+\half) \tilde{s}_m \<\<\sigma_{m+k}\>\> +
\tfrac\hbar8 \sum_{i+j=k-1} \Bigl( \<\<\sigma_i\sigma_j\>\> +
\<\<\sigma_i\>\> \<\<\sigma_j\>\> \Bigr) .
\end{align*}
We now show, following Dijkgraaf et al., that these constraints are a
formal consequence of Witten's conjecture and the puncture equation
$z_{-1}=0$. The proof of Theorem \ref{DVV} does not use the puncture
equation, and holds for any solution of the KdV hierarchy.
\begin{theorem} \label{DVV}
The recursion \eqref{f} implies the recursion
$$
\p \HH \p z_{k-1} = \p^3 z_k , \quad k\ge0 .
$$
\end{theorem}

The constraints $z_k=0$ follow from Theorem \ref{DVV} by induction from
this recursion, starting with the puncture equation $z_{-1}=0$:
\begin{enumerate}
\item the induction hypothesis $z_{k-1}=0$ and Theorem \ref{DVV} imply that
$\p^3z_k=0$;
\item applying Theorem \ref{basic}, we conclude that $z_k=0$.
\end{enumerate}

In the presence of the puncture equation, both Witten's conjecture and the
Virasoro conjecture determine the Gromov-Witten potential uniquely; we
conclude that these conjectures are equivalent.

There are a number of other proofs that Witten's conjecture
\eqref{witten-conjecture} implies the Virasoro constraints: Goeree
\cite{Goeree} and Kac and Schwartz \cite{KS} give proofs using vertex
operators, while La \cite{La} uses the theory of Lie-B\"acklund
transformations. However, we have chosen to present the original proof,
since it is completely elementary.

A direct proof that Kontsevich's integral representation of the potential
function $Z$ satisfies the Virasoro constraints was given by Witten
\cite{W}. Later, simpler derivations were given by Gross and Newman
\cite{GN} and by Itzykson and Zuber \cite{IZ}.

\raggedbottom

\subsection{Proof of Theorem \ref{DVV}}
We leave the proof that $\HH\p z_{-1}=\p^2z_0$ to the reader. Turning to
$k>0$, we divide the calculation of $\p^2z_k$ into three parts:
\begin{align*}
\text{I} &= \p^2 \sum_{m=0}^\infty (m+\half) \tilde{s}_m
\<\<\sigma_{m+k}\>\> = \sum_{m=0}^\infty (m+\half) \tilde{s}_m
\p^2 \, \<\<\sigma_{m+k}\>\> + \p\<\<\sigma_k\>\> \\
&= \sum_{m=0}^\infty (m+\half) \tilde{s}_m \, \HH\p
\<\<\sigma_{m+k-1}\>\> + \p\<\<\sigma_k\>\> \\
&= \HH\p \sum_{m=0}^\infty (m+\half) \tilde{s}_m
\<\<\sigma_{m+k-1}\>\> + \p\<\<\sigma_k\>\> - \bigl( \tfrac{\hbar}{4} \p^3
+ \u\p + \tfrac{1}{4} \u' \bigr) \<\<\sigma_{k-1}\>\> ;
\end{align*}
\begin{align*}
\text{II} &= \p^2 \sum_{i+j=k-1} \<\<\sigma_i\sigma_j\>\> =
\sum_{i+j=k-2} \p_j \p^2 \<\<\sigma_{i+1}\>\> + \p^3\<\<\sigma_{k-1}\>\> \\
&= \sum_{i+j=k-2} \p_j \HH \p \<\<\sigma_i\>\> +
\p^3\<\<\sigma_{k-1}\>\> \\ {} &= \HH\p \sum_{i+j=k-2}
\<\<\sigma_i\sigma_j\>\> + \p^3\<\<\sigma_{k-1}\>\> \\
& \quad + \sum_{i+j=k-2} \Bigl( \tfrac\hbar2 \p\<\<\sigma_i\>\>
\p^3\<\<\sigma_j\>\> + \hbar\, \p^2\<\<\sigma_i\>\> \p^2\<\<\sigma_j\>\>
\Bigr) ; \\
\text{III} &= \p^2 \sum_{i+j=k-1} \<\<\sigma_i\>\> \<\<\sigma_j\>\>
= 2 \sum_{i+j=k-1} \<\<\sigma_i\>\> \p^2 \<\<\sigma_j\>\> + 2
\sum_{i+j=k-1} \p \<\<\sigma_i\>\> \p \<\<\sigma_j\>\> \\
&= 2 \sum_{i+j=k-2} \<\<\sigma_i\>\> \p^2 \<\<\sigma_{j+1}\>\> +
\tfrac{2}{\hbar} \u' \<\<\sigma_{k-1}\>\> + 2 \sum_{i+j=k-1} \p
\<\<\sigma_i\>\> \p \<\<\sigma_j\>\> \\ {} &= 2 \sum_{i+j=k-2}
\<\<\sigma_i\>\> \HH\p \<\<\sigma_j\>\> + \tfrac{2}{\hbar} \u'
\<\<\sigma_{k-1}\>\> + 2 \sum_{i+j=k-1} \p \<\<\sigma_i\>\> \p
\<\<\sigma_j\>\> \\ {} &= \HH\p \biggl( \sum_{i+j=k-2} \<\<\sigma_i\>\>
\<\<\sigma_j\>\> \biggr) + \tfrac{2}{\hbar} \u' \<\<\sigma_{k-1}\>\> + 2
\sum_{i+j=k-1} \p \<\<\sigma_i\>\> \p \<\<\sigma_j\>\> \\ {} &\quad -
\sum_{i+j=k-2} \Bigl( \hbar \, \p\<\<\sigma_i\>\> \p^3\<\<\sigma_j\>\> +
\tfrac{3\hbar}{4} \p^2\<\<\sigma_i\>\> \p^2\<\<\sigma_j\>\> +
2\u \, \p\<\<\sigma_i\>\> \p\<\<\sigma_j\>\> \Bigr) .
\end{align*}
Combining these calculations, we see that
$$
\p^2 z_k - \HH\p z_{k-1} = \text{I} + \tfrac{\hbar}{8} ( \text{II} +
\text{III} ) - \HH\p z_{k-1} = a + b ,
$$
where
\begin{align*}
a &= \p\<\<\sigma_k\>\> - \bigl( \tfrac{\hbar}{8} \p^3 + \u\p \bigr)
\<\<\sigma_{k-1}\>\> + \tfrac{\hbar}{4} \sum_{i+j=k-1} \p \<\<\sigma_i\>\>
\p \<\<\sigma_j\>\> \\
b &= - \tfrac{\hbar}{2} \sum_{i+j=k-2} \Bigl( \tfrac{\hbar}{8} \p
\<\<\sigma_i\>\> \p^3 \<\<\sigma_j\>\> - \tfrac{\hbar}{16} \p^2
\<\<\sigma_i\>\> \p^2  \<\<\sigma_j\>\> + \half \u \p \<\<\sigma_i\>\> \p
\<\<\sigma_j\>\> \Bigr) .
\end{align*}
It follows that $\p\bigl(\p^2z_k - \HH\p \bigr)z_{k-1} = \p a + \p b$,
where
\begin{align*}
\p a &= \p^2\<\<\sigma_k\>\> - \bigl( \tfrac{\hbar}{8} \p^3 + \u\p + \u'
\bigr) \p\<\<\sigma_{k-1}\>\> + \tfrac{\hbar}{2} \sum_{i+j=k-1} \p
\<\<\sigma_i\>\> \p^2 \<\<\sigma_j\>\>  \\ 
&= \tfrac{\hbar}{2} \sum_{i+j=k-2} \p \<\<\sigma_i\>\> \p^2
\<\<\sigma_{j+1}\>\> \\
\p b &= - \tfrac{\hbar}{2} \sum_{i+j=k-2} \p \<\<\sigma_i\>\> \HH\p
\<\<\sigma_j\>\> .
\end{align*}
Thus $\p(a+b)=0$, completing the proof of Theorem \ref{DVV}.
\qed

\flushbottom

\bigskip

Dubrovin and Zhang \cite{DZ} have proved an analogue of Theorem \ref{DVV}
for any Frobenius manifold $\phase$, but only in genus $0$ --- in
particular, for the Gromov-Witten invariants of any smooth projective
variety $V$. The recursion operator of the KdV hierarchy has the genus $0$
limit
$$
\HH_0 = \lim_{\hbar\to0} \HH = u\p + \half u' ,
$$
and its analogue in the general case is given by the formula
$\DC\p+(\mu+\half)\p\CC$.  Like the operator $u\p+\half u'$, this operator
is Hamiltonian: that is, the bilinear differential operator
\newcommand{\lp}{\{\!\{}
\newcommand{\rp}{\}\!\}}
$$
\lp u_a(x) , u_b(y) \rp_0 = \DC_{ab} \delta'(x-y) + (\mu_b+\half)
\p\CC_{ab} \delta(x-u)
$$
on the algebra of differential polynomials $\Q\{u_a\}$ is a Poisson
bracket. Furthermore, together with the Poisson bracket
$$
\{ u_a(x) , u_b(y) \}_0 = \delta'(x-y) ,
$$
it generates a pencil of Poisson strctures.

Dubrovin and Zhang have conjectured that, if the Frobenius manifold
$\phase(V)$ is semisimple, these Poisson brackets are the genus $0$ limits
of Poisson brackets of the form
\begin{align*}
\lp u_a(x) , u_b(y) \rp &= \lp u_a(x) , u_b(y) \rp_0 + \sum_{g=1}^\infty
\sum_{i=0}^{2g+1} \hbar^g k_{i,g} \delta^{(i)}(x-y) , \\
\{ u_a(x) , u_b(y) \} &= \{ u_a(x) , u_b(y) \}_0 + \sum_{g=1}^\infty
\sum_{i=0}^{2g+1} \hbar^g h_{i,g} \delta^{(i)}(x-y) .
\end{align*}
where $k_{i,g},h_{i,g} \in \Q\{u^a\} \o \End(H(V))$, that these two Poisson
brackets generate a pencil of Poisson structures, that the hierarchy of
commuting Hamiltonian flows associated to the functions
$\<\<\tau_{0,a}\tau_{0,0}\>\>$ is $\<\<\tau_{n,a}\tau_{0,0}\>\>$, $n\ge0$,
 and that the Virasoro constraints define Lie-B\"acklund transformations of
this hierarchy; they have verified this conjecture up to genus $1$. This
makes it plausible that, when the Frobenius manifold $\phase(V)$ is
semisimple, the Gromov-Witten invariants in \emph{all} genera are
determined by those in genus $0$ together with the Virasoro constraints.

An approach to calculating the higher genus Gromov-Witten invariants has
been outlined by Eguchi and Xiong \cite{EX}, and worked out for
$\mathbb{P}^2$ in genus $2$ and $3$; they combine the Virasoro constraints
with equation which follow, by the topological recursion relations of
\cite{genus2,KM1}, from the obvious fact that any monomial
$\psi_1^{k_1}\dots\psi_n^{k_n}$ on $\Mbar_{g,n}$ vanishes if
$k_1+\dots+k_n>3g-3+n$.

\section{The Virasoro conjecture for Calabi-Yau manifolds} \label{CY}

A Calabi-Yau manifold $V$ is a smooth projective variety such that
$c_1(V)=0$ and $H^1(V,\C)=0$. In this section, we prove the Virasoro
conjecture for Calabi-Yau varieties; this generalizes unpublished results
of S. Katz for threefolds, while the suggestion to consider other
dimensions was made by J. Bryan. These instances of the Virasoro conjecture
are in fact a little dull: they impose no constraints on the Gromov-Witten
invariants of $V$.
\begin{theorem}
If $V$ is a Calabi-Yau variety, the Virasoro constraints $z_{k,g}=0$ hold.
\end{theorem}
\begin{proof}
If $V$ is a holomorphic symplectic manifold with $h^{2,0}=1$ (this includes
K3 surfaces, as well as abelian surfaces), the Gromov-Witten invariants of
$V$ vanish except possibly in degree $\beta=0$ (Behrend and Fantechi
\cite{BF}), while the Virasoro conjecture holds in degree $\beta=0$ by the
explicit calculations of \cite{GP}. Thus, we may assume that $r\ge3$.

If $c_1(V)=0$, the formula for the virtual dimension of $\Mbar_g(V,\beta)$
is very simple:
$$
\vdim \Mbar_g(V,\beta) = (3-r)(g-1) .
$$
Fixing $g>0$ and starting with Hori's equation $z_{0,g}=0$, we will prove
that $z_{k,g}$ vanishes by induction on $k$, using the fact that
$\vdim\Mbar_g(V,\beta)\le0$.

By Lemma \ref{downward}, $\CL_{-1}z_{k,g}=-(k+1)z_{k-1,g}$, which vanishes
under the induction hypothesis $z_{k-1,g}=0$. Writing this equation out
explicitly, we see that
$$
\p_{0,0} z_{k,g} = \sum_{m=0}^\infty t^a_{m+1} \p_{m,a} z_{k,g} .
$$
This shows that $z_{k,g}$ is determined by its restriction to the subscheme
$\{t_0^0=0\}$ of $\PHASE(V)$. Let $i$ be the embedding
$\{t_0^0=0\}\hookrightarrow\PHASE(V)$.

From the dimension equation
$$
\sum_{m=0}^\infty ( p_a + m - 1 ) t_m^a \p_{m,a} \<\<~\>\>_g =
\vdim\Mbar_g(V,\beta) \* \<\<~\>\>_g ,
$$
we see that
$$
\sum_{m=0}^\infty ( p_a + m - 1 \bigr) t_m^a \p_{m,a} z_{k,g} = (
\vdim\Mbar_g(V,\beta) - k ) z_{k,g} .
$$
This equation has an anti-holomorphic partner, obtained by replacing $p_a$
by $q_a$:
$$
\sum_{m=0}^\infty ( q_a + m - 1 \bigr) t_m^a \p_{m,a} z_{k,g} = (
\vdim\Mbar_g(V,\beta) - k ) z_{k,g}
.
$$
The vector fields entering into these equations are tangential to the
subscheme $\{t_0^0=0\}$; adding them together, we conclude that
$$
\sum_{m=0}^\infty ( \half p_a + \half q_a + m - 1 \bigr) t_m^a \p_{m,a}
i^*z_{k,g} = ( \vdim\Mbar_g(V,\beta) - k ) i^*z_{k,g} .
$$
The left-hand side of this equation is positive semi-definite in the
monomial basis of the algebra of functions on $\{t_0^0=t_1^0=0\}$, since
$H^1(V,\C)=0$; we conclude that $z_{k,g}$ vanishes.
\end{proof}

\end{document}